\input amstex
\documentstyle{amsppt}

 \voffset 1.35 truein
 \hoffset 0.845 truein
 \input defs
 \topmatter
\title%\titlebf
\centerline{Isospectral Pairs of Metrics
 on Balls, Spheres,}
\centerline{
 and other Manifolds}
\centerline{  
with Different Local Geometries} 
\endtitle
\author\tf Z.I. Szab\'o\endauthor
\thanks Research partially supported by 
the CUNY grants 9-92720 and 9-91909.
\endthanks
\headLine{Z. I. SZAB\'O}{ISOSPECTRAL PAIRS ON BALLS AND SPHERES}
\address
Zolt\'an Imre Szab\'o
\endaddress
\address
\noindent
CUNY,
Lehman College,
Bronx, NY 10468-1589,
USA
\endaddress
\address
R\'enyi Alfr\'ed Institute of Math.,
P. O. Box 127,
H-1053 Budapest,
Hungary
\endaddress
\address
Email: {\rm zoltan\@alpha.lehman.cuny.edu}
\endaddress

\date
January 2000
\enddate

\abstract The first isospectral pairs of metrics 
are constructed
on the most simple simply connected domains, namely, on balls
and spheres. 
This long standing problem, concerning the existence of such
pairs, has been solved by a 
new method called "Anticommutator Technique".
 Among the wide range of such pairs, the most Striking  
Examples 
are provided on the spheres 
$S^{4k-1}$, where $k\geq 3\,$.
One of these metrics 
 is homogeneous 
(since it is the metric on the geodesic sphere 
of a 2-point homogeneous
space),
while the other is locally
inhomogeneous. These examples demonstrate the surprising fact that no
information
about the isometries  
 is encoded in the spectrum of Laplacian acting on
functions.
In other words,
"The group of isometries, even the local homogeneity property,
 is lost to the "Non-Audible" in the debate of
"Audible versus Non-Audible Geometry"."
\endabstract

\endtopmatter

 \rightline{\it To the memory of my son D\'aniel (1976-1999)} 
\vskip 0.85 cm

Research in Spectral Geometry has been booming in the
last two decades.
Since Milnor's
\cite{M,1964} 
 first example of a pair
of 16 dimensional isospectral but non-isometric flat tori, 
many new examples
were constructed.
See some of the most important
results in 
\cite{V,1980}; \cite{I,1980};
\cite{GW,1984,1986}; \cite{S,1985}; \cite{Bu,1986};
\cite{BT,1987}
\cite{DG,1989}; \cite{GWW,1992}, etc.
However,
all of these previously 
constructed isospectral pairs consist of locally isometric
spaces and they differ from each other just in global shape;
global isometries do not exist in these cases but local
isometries do.

The first examples of isospectral metrics with different
local geometries were constructed  
by the author \cite{Sz3,1992} and by C. Gordon \cite{G1,1993}.
In these new examples Gordon covered the closed case 
while the author constructed metrics on the product of a ball
by a torus (the boundary of this manifold is the product of a
sphere by the same torus) such that the metrics were both
Dirichlet and Neumann isospectral. Yet the spaces were 
locally non-isometric.
\par To the construction of these examples we used special
2-step nilpotent Lie algebras. 
In the following we briefly introduce some of those
 basic notions 
what are really needed 
for a clear exposition.

\begingroup\eightpoint
A general metric
2-step nilpotent Lie algebra is described
by the system 
$$
\big\{\bn=\bv\oplus\bz,\<,\>,\bJ_Z\big\},
\tag 0.1
$$
where $\bz=[\bn,\bn]$ is the center of the Lie algebra, $\<,\>$
is an inner product on $\bn$ such that $\bv\bot\bz$ furthermore
the one-to-one mapping $\bJ : \bz \to \End(\bv)$ is defined by
the formula
$$
\<\bJ_Z (X),Y\>=\<[X,Y],Z\>\quad ,\quad \forall X,X'\in \bv \quad ,\quad
\forall Z\in\bz.
\tag 0.2
$$

The space $\bv$ resp. $\bz$ is called also
X-space resp. $Z$-space. The linear space of the
endomorphisms $\bold J_Z,\forall Z\in\bold z$ is denoted by 
$\bold {J_z}$.
Special 2-step nilpotent Lie algebras are the so called Heisenberg-type
algebras,
characterized by the property 
$\bold J_Z^2=-|Z|^2id,\forall Z\in\bold z$ \cite{K}.

Any 2-step nilpotent Lie group $N$ can be extended into a solvable group
$SN$ defined on the half space $\bold n\times\bR_+$. The corresponding
formulas of this extension are described in (3.16) and (3.17). 
The natural invariant metrics
 on the solvable extensions of the Heisenberg-type groups
 are harmonic \cite {DR}. One can find
more about harmonic manifolds in \cite {Sz1,2,3}.

In case of the special 2-step nilpotent Lie 
algebra $\bold H^{(a,b)}_3$,
 used for constructing the above new examples,
the $Z$-space is $\bR^3$ (considered as the space
of imaginary quaternions) and the $X$-space is the $4(a+b)$-
dimensional  quaternionic space $\bH^{a+b}$,
where $\bH=\bR^4$ is the space of quaternions and 
$a,b\in\bold N$ are natural
numbers.
The inner product on $\bn$ is the natural inner product of the
space $\bn=\bR^{4(a+b)}\times\bR^3$. The endomorphisms
$\bJ_Z, Z\in\bR^3$
 are introduced by 
using both left and right products as follows:
$$
\bJ^{(a,b)}_Z(X_1,\dots,X_a,X_{a+1},\dots,X_{a+b})
=\big(X_1Z,\dots,X_aZ,ZX_{a+1},\dots,ZX_{a+b}\big).
\tag 0.3
$$
The corresponding Lie algebras and Lie groups are denoted by
$\bold h^{(a,b)}_3$ and \ $\bH^{(a,b)}_3$.
The metric defined by the left-invariant extension of the inner
product $\<,\>$ of the tangent space
$\bold h^{(a,b)}_3=T_0(H^ {(a,b)}_3)$ onto the group
is denoted by $g^{(a,b)}$. 

The spectral investigations will be performed on four type of domains.
These are the Ball$\times$Torus-type domains, the pure Ball-type
domains, and also their Sphere$\times$Torus- resp. the pure Sphere-type
boundaries of these domains will be considered. 
They are described as follows.

Let $\Ga_Z$ be a full lattice (discrete subgroup) on the Z-space
and let $B_\de$ be the euclidean ball of radius $\de$ around 
the origin of the X-space. The factor manifold $\Ga_Z\setminus
\bH^{(a,b)}_3$ can be considered as a principal fiber bundle 
such that the
fibers are naturally identified with the torus $T^3=\Gamma /\bR^3$ and
the base space is naturally identified with the X-space. 
When the fibers are
considered only over $B_\de$, the bundle is denoted by
$(B_\de,T^3)$. The boundary of this manifold is the principal 
fiber bundle
 $(S_\de,T^3)$, where
 the sphere $S_\de =\pa B_\de$ is the base space of this bundle. This
construction can be established on any 2-step nilpotent group.

Besides the above considered Ball$\times$Torus- and 
Sphere$\times$Torus-cases
 we consider also such compact domains whose boundaries are
 described
by equations of the form $f\big(|X|,|Z|\big)=0$. In these cases the
 domains
are diffeomorphic to euclidean  balls and the boundaries are 
diffeomorphic to euclidean spheres. We call
 these
cases {\it Ball-cases\/} and {\it Sphere-cases\/}. 

The corresponding domains and their boundaries on the solvable
extensions are introduced by a natural modification. 
\endgroup

Then we have
\proclaim{Theorem \cite{Sz3}} 
The metrics $g^{(a,b)}$
with the same $(a+b)$ but with different $ab$ are locally
non-isometric, yet they are both Dirichlet and Neumann isospectral
on the manifolds $(B_\de,T^3)$. 

Beside proving the isospectral
property by constructing an intertwining operator,
also an explicit computation of the spectrum 
is provided in the paper.
\endproclaim
 
Gordon and Wilson generalized this construction \cite{GW3} to
obtain continuous families of isospectral metrics such that 
the family members have different local geometries (namely the scalar
curvature is different on these spaces).  
Their Main Theorem states
\proclaim{Theorem \cite{GW3}}
Let $\big\{\bold n=\bv\oplus\bz,\<,\>,\bJ_Z\big\}$ and 
$\big\{\bold {n'}=\bv\oplus\bz,\<,\>,\bold {J'}_Z\big\}$
be Lie algebras 
defined over the same space $\bv\oplus\bz$ such that $\bJ_Z$ and 
$\bold {J'}_Z$ are isospectral (conjugate), $\forall Z\in\bz$. Then the 
left-invariant metrics $g$ and $g'$ are both Dirichlet and Neumann 
isospectral on any of the Ball $\times$ Torus-type manifolds
described 
above.
\endproclaim

Notice that the endomorphisms introduced in (0.3) are spectrally
equivalent in the sense of the above theorem on spaces with the 
same $a+b$. It should be mentioned that 
 the
existence of non-trivial isospectral deformations is proven
in \cite{GW3} only on spaces where this torus is 2-dimensional. It 
is still an open question if this construction provides 
non-trivial continuous 
isospectral deformations on spaces having higher 
dimensional Z-spaces. This points out one of the differences between 
the constructions of \cite{GW3} and \cite{Sz3}. It should be emphasized
too that 
 this generalization strictly concerns the Ball$\times$Torus-type 
manifolds and it is not valid, for instance, on pure Ball-type domains.

 Further developed  this field by Gordon, Gornet, Schueth, Webb
and Wilson \cite{GGSWW} by noticing that the metrics $g_t$
restricted to the boundary $S_\de\times T^2$ of the spaces
constructed in \cite{GW3} are isospectral with different local
geometries. The new feature of these examples is that these
spaces are closed and locally inhomogeneous.

Independently also this author noticed that {\it the metrics 
$g^{(a,b)}$ restricted to the boundary $S_\de\times T^3$
of his examples are isospectral with different local 
geometries.} (This observation is included in a revised version
of the article \cite{Sz3}.) The new feature of these
examples is that one of these spaces (belonging to $g^{(0,a+b)}$)
is homogeneous while the other spaces (belonging to $g^{(a,b)}$,
where $ab\neq 0$) are locally inhomogeneous. This is the first 
demonstration of the fact that {\it one can not hear the local
homogeneity property of closed Riemannian manifolds.}

 Recently D. Schueth \cite{Sch} constructed continuous 
families of isospectral Riemannian metrics on simply connected
manifolds such that the different members of a family have different
scalar curvature. In her construction she enlarges the torus
$T^2$ of the spaces considered in \cite{GGSWW} into a compact
simply connected Lie group $G$ such that $T^2$ is a maximal torus
of $G$. Then she obtains her examples by an extension of the metrics
$g_t$ of $S_\de\times T^2$ onto the enlarged manifold $S_\de\times G$.
In special cases she gets examples on the product of spheres. The space
with the lowest dimension is constructed on
$S^4\times S^3\times S^3$. These are the first examples for 
isospectral non-isometric metrics constructed on closed simply 
connected manifolds.
 
  Most recently Gordon and this author
\cite{GS} constructed continuous families of negatively curved
manifolds with boundaries such that the metrics of the distinct family
 members are 
isospectral, yet
having different local geometries. This result contrasts the Rigidity
Theorem of Crocke and Sharafutdinov\cite{CS}, asserting that such 
continuous
deformations 
do not exist on closed manifolds.

It should be mentioned
that all the above constructions involve only the spectrum of the 
Laplacian acting on functions. The manifolds used for 
these constructions
are exclusively the manifolds described in
 the Ball$\times$Torus-cases or in 
the Sphere$\times$Torus-cases. No examples of non-trivial
isospectral pairs are constructed 
 on the most simple simply connected
manifolds, namely, on balls or on spheres so far. 

Examples concerning the 
spectra of forms
or spectrum of the Dirac operator are also completely missing from
the literature. The most far reaching statement from this point 
of view is that the metrics constructed by Schueth are not isospectral 
on 1-forms \cite{Sch}. 

We extend our investigations in these directions. The main goal of 
this paper is to find the first non-trivial isospectral pairs of
spaces constructed on the most simple domains, namely, on balls and
spheres. However, we scrutinize also the question about
the importance of using spectrally equivalent endomorphism spaces
for constructing isospectral metrics. It turns out that this 
property is
important in case of Ball$\times$Torus-type domains and on their
boundaries but it
is not important, at all, on the Ball-type domains and on their
Sphere-type boundaries. 
A brief 
outline of these considerations is as follows.

First discrete-type, so called {\it $\si$-deformations} are introduced
 on 2-step 
nilpotent Lie 
algebras.
These deformations are generalizations of the special 
deformations used for constructing the non-isomorphic Lie algebras
$
\bold h^{(a,b)}_3 $ 
in (0.3). Also these generalized deformations provide
 spectrally equivalent endomorphism metrics. 
 The more precise details are as follows.

Consider a Lie algebra $\big\{\bn=\bv\oplus\bz,\<,\>,\bJ_Z\big\}$ 
and let
$\si$ be an involutive orthogonal transformation on $\bv$
commuting with the endomorphisms $\bJ_Z,\forall Z\in\bz$. 
Then the $\si$-deformation of this
Lie algebra is defined by
 $\big\{\bn_\si =\bv\oplus\bz,\<,\>,\bJ^{\si}_Z\big\}$,
 where 
only the endomorphisms $\bJ^{\si}_Z$ are different in the 
new system and they
are 
defined by 
$\bJ^\si _Z=\si\bJ_Z,\forall Z\in\bz$.
The corresponding 
left invariant metrics on the $\sigma$-deformed groups 
are denoted by $g^\si$.

More general deformations are the so called {\it Partial 
$\si$-deformations} (they are called also 
{\it $P\si$-deformations}). In this case we perform 
$\si$-deformations only on
the elements of a subspace 
$\bJ_{\bold s}\subset\bJ_\bz$,
while keeping the orthogonal 
complement $\bJ_{\bold s^\perp}$ 
unchanged. 

Though the Partial $\sigma$-deformations
 form a much wider class than the $\si$-deformations, in very particular
cases a $\si$-deformation can be
 equivalent to a $P\si$-deformation. This
is the case
 when the endomorphism space contains at least
a 1-dimensional subspace,
$\bJ_{\bold A}$, 
of the so called anticommutators. A non-degenerated
endomorphism $A\in \bold {J_z}$ is an anticommutator
iff $A\circ B=-B\circ A,\forall B\in \bold J_{A^\perp}$.
 A Partial $\si$-deformation established by the deformation $A^\si$ 
 of an anticommutator $A$,
is called $\sigma_A$-deformation.
In Reduction Theorem 4.1 we prove that a
$\si$-deformation is always equivalent to the $\sigma_A$-deformation,
 if the endomorphism space contains an anticommutator $A$.  

 The general $\si$-deformations provide isospectral endomorphisms
in the sense of Theorem \cite{GW3}, therefore 
{\it they provide a wide range
of new examples for non-trivial isospectral pairs on the 
Ball$\times$Torus-type
domains.} However, the general $sigma$-deformations do not provide
isospectral metrics on the pure Ball-type domains and on their 
Sphere-type boundaries. Our main goal in this paper is to find
"such/those" $sigma$-deformations which lead to isospectral metrics
also on these simple domains. It turns out that the particular 
$\si_A$-deformations are "such" deformations. More precisely, we have
\proclaim{Main Theorem 4.3} Let 
$\bold{J_z}=\bold A\oplus\bold A^\perp$ 
be an Endomorphism Space
with the anticommutator $A$ and let $\si$ be an involutive orthogonal
transformation on the X-space commuting with all of the endomorphisms
from $\bold{J_z}$. Then the metrics $g$ and $g^\si$, determined by
the Endomorphism Spaces $\bold{J_z}$ and 
$\bold{J^\si_z}\simeq\bold A^\si\oplus\bold A^\perp$, are both Dirichlet
and Neumann isospectral on the Ball-type domains of the considered
manifolds.

This Isospectrality Theorem holds also on the Sphere-type boundaries
of these domains and both versions of this Theorem can be established
also on the solvable extensions of 2-step nilpotent groups.
\endproclaim 

The Theorem is proven by constructing intertwining operators. It says
that,
on balls and spheres, 
isospectral metrics can be constructed
by appropriate deformations of an anticommutator in an Endomorphism
Space. 
In this paper we consider mostly the $\sigma_A$-deformations of
an anticommutator. 
However, the proofs given in this case make it
apparent that they are
 "working" in a much more general situation, namely,
when the anticommutator is deformed only isospectrally. The full
 exploration of
all options for such deformations  
will be completed in the continuation \cite{Sz4} of this paper. There
we discover also a cornucopia of the so called $SpeCh_A$-deformations which
provide spectrally inequivalent Endomorphism
 Spaces and therefore a completely
new type of examples. These $SpeCh_A$-deformations show
that our method is drastically
forking from all of the methods applied on this field so far.
 This completely new method is called {\it constructions 
of isospectral metrics by means
of isospectral deformations of an anticommutator in an Endomorphism
Space.} Or, we simply call it {\it Anticommutator Technique}.
This new method inherits 
no "tricks" from the previous investigations performed on 
this field, since all the tools applied
on the Ball$\times$Torus-type manifolds so far break down on the Ball-type
domains and on their Sphere-type boundaries.

We provide also 
the corresponding {\bf Non-Isometry Proofs}
 in this paper.
Thus a wide range of isospectral pairs of metrics are constructed also
on balls and spheres such that the metrics in a pair
 define different local geometries. By any means,
these are the very first examples of isospectral pairs of metrics
constructed on these most simple simply connected domains.

The most {\bf Striking Examples} are constructed on the geodesic
spheres of the Harmonic Spaces 
$S\bold H^{(a,b)}_3$). 
Since the geodesic spheres
are homogeneous
on the 2-point homogeneous space $S\bold H^{(a+b,0)}_3$ 
and the geodesic spheres 
 are inhomogeneous
on the other considered harmonic spaces,
these examples provide 
isospectral pairs of metrics constructed on a sphere such that one of them
is a homogeneous and the other is a locally inhomogeneous metric.
These examples demonstrate the fact that the group of isometries is
lost to the "Non-Audible" in the debate of "Audible versus Non-Audible
Geometry".

The question
 about the isospectrality of the geodesic balls
and spheres on the above harmonic spaces was
 asked in \cite {Sz2}. The investigations on this field took several
unexpected turns, till the Isospectrality Theorem has been 
established in the 
general form stated in this paper.

Finally,
we would like to mention two related Statements. 
Each of them is
suggestive concerning the possible form of the
Isospectrality
Theorem 
considered on forms.
\proclaim{Propositions 5.3-5.4}
 Two $\sigma$-equivalent metrics have 
identical Ricci curvature at the origin. However, 
they have only isotonal
Riemannian curvature. 
\endproclaim

(Two operators are said to be {\it isotonal\/} 
if the elements of the spectra
 are the 
same, but the multiplicities may be different. 
(Of course, 0-multiplicity is
 not allowed 
in this counting.) 
We use this notion 
with respect to the Laplacian as well as to the
 curvature operator ${R_{ab}}^{cd}$ acting as a
symmetric endomorphism on the space of 2-forms.)

This statement contrasts 
the following statement:
\proclaim{Proposition 5.1-5.2}
 Two $\si$-equivalent metrics are Laplace-isospectral
on the invariant 1- and 2-forms.
\endproclaim
(Formula (1.18) shows the involvement of the curvature operator into
the Laplacian acting on forms, showing also the 
contrasting nature of the above two statements.)

This is the opening paper of related articles.
Since we did not want to cram to much technical details in one paper,
we divided these technicalities between this paper
and article \cite{Sz4}. Yet, we tried to make them 
 as independent as it
was possible. In this paper the proofs are complete
on manifolds with boundaries (i.e. on Ball$\times$Torus-type and on
 Ball-type manifolds) and some of the technical details, what are 
 needed on the boundary manifolds, are
discussed in \cite{Sz4}. In some extend this remark concerns also the
Striking Examples, where the Non-Isometry Proof is based on an explicit
computation of the eigenvalues of the Ricci operator. These computation
is incorporated into \cite{Sz4}. All the other 
details are self-contained. 

I would like to express my gratitude to Carolyn Gordon for the helpful
 discussions resulting in a clearer exposition 
and to Attila M\'at\'e for helping me
in the preparation of this article.

  \head
 {\titlebf 1. \S\quad Two-step nilpotent Lie algebras}
\endhead

In the Introduction we reviewed some of the basic materials about
2-step nilpotent Lie algebras. They are described by a system 
$$
\big\{
\bold n=\bv\oplus\bz,\<,\>,\bJ_Z\big\},
$$
 as it is explained in
(0.1) and
(0.2). A great deal about these Lie algebras can be found in \cite{E}. 
Such 
a Lie algebra is well defined by the endomorphisms $\bJ_Z$ and by the 
inner product $\<,\>$ defined on $\bold n$. The linear space
of these endomorphisms is denoted by $\bJ_\bz$ and for an X-vector 
$X\in\bold v$,
$\bJ_\bz(X)$ denotes the linear subspace of $\bold v$ 
spanned by the vectors of the form
$\bJ_Z(X),\forall Z\in\bz$.

The exponential map identifies the Lie group $N$ with the space 
$\bv\oplus\bz$ and the group product on $N$ is described by the formula
$$
(X,Z)\cdot (X',Z')=\big(X+X', Z+Z'+1/2[X,X']\big).
\tag 1.1 $$ 

Consider the orthonormal bases $\big\{E_1;\dots;E_k\}$ and 
$\big\{e_1;\dots;e_l\big\}$ on the X- and on the Z-space respectively.
The corresponding coordinate systems defined by these bases are 
denoted by
$\big\{x^1;\dots;x^k\big\}$ and $\big\{z^1;\dots;z^l\big\}$. We
extend the vectors $E_i;e_{\alpha}$ into the left-invariant
vector fields
$$
\gathered
\bold
X_i=\pa_i + \frac 1 {2} \sum_{\alpha =1}^l
\<[X,E_i],e_{\alpha} \>  \pa_\alpha = 
\\
=\pa_i + \frac 1 {2} \sum_{{\alpha} =1}^l \<
\bold J_\alpha\big(X\big),E_i\>
\pa_\alpha\quad ; \quad
\bold Z_{\alpha}=\pa_\alpha,
\endgathered
\tag 1.2
$$
where $\pa_i=\pa /\pa x^i$, $\pa_\alpha=\pa/\pa z^\alpha$
and $\bold J_\alpha =\bold J_{e_\alpha}$.
\par By using these
vector fields, we establish some explicit formulas on
2-step nilpotent Lie groups.

The covariant
derivative can be computed by the well known formula
$$
\<\na _PQ,R\>=\frac 1 {2} \big\{\<P,[R,Q]\>+\<Q,[R,P]\>+
\<[P,Q],R\>\big\},
\tag 1.3
$$
where $P,Q,R$ are considered as elements of the Lie algebra $\bold n$
(or, they can be considered as  invariant vector fields on $N$). 
By this formula 
we get
$$
\na_XX^*=\frac 1 {2} [X,X^*];\na_XZ=\na_ZX=-\frac 1 {2}
\bold J_Z\big (X\big);\na_ZZ^*=0.
\tag 1.4
$$

 The Laplacian $\Delta$ acting on functions can be 
explicitly computed by using
the formula
$$
\De=\sum_{i=1}^k\big(\bold X_i^2-\na_{\bold X_i}\bold X_i\big)
+\sum_{\al=1}^l\big (\bold Z_{\al}^2-\na_{\bold Z_{\al}} 
\bold Z_{\al}\big ).
\tag 1.5
$$
 Then one obtains the following explicit formula
$$
\De=\De_X+\De_Z+\frac 1 {4} \sum_{\al,\beta =1}^l \<\bold J_\al
\big (X\big),\bold J_\beta\big (X\big)\>
 \pa_{\alpha\beta}^2
+\sum_{\al=1}^l\pa_\alpha D_\al \bullet,
\tag 1.6
$$
where $D_\al \bullet$ means differentiation (directional derivative)
with respect to the vector field
$$
D_\al : X \to \bold J_{\al}\big (X\big )
\tag 1.7
$$
tangent to the X-space, furthermore 
$\pa_{\alpha\beta}=\pa^2/\pa z^\alpha
\pa z^\beta$.\par We compute also some other basic objects explicitly.
The Riemannian curvature tensor on 2-step 
nilpotent groups are explicitly
computed in \cite {E}. By those formulas we have:
$$
\gathered
R(X,Y)X^*=\frac 1 {2} \bold J_{[X,Y]}(X^*) -
 \frac 1 {4} \bold J_{[Y,X^*]}
(X) + \frac 1 {4} \bold J_{[X,X^*]}(Y); \\
R(X,Y)Z=-\frac 1 {4} [X,\bold J_Z(Y)]+\frac 1 {4} 
[Y,\bold J_Z(X)]\quad ;
\quad R(Z_1,Z_2)Z_3=0; \\
R(X,Z)Y=-\frac 1 {4} [X,\bold J_Z(Y)]\quad ; 
\quad R(X,Z)Z^*=-\frac 1 {4}
\bold J_Z\bold J_{Z^*}(X); \\
R(Z,Z^*)X=-\frac 1 {4} \bold J_{Z^*}\bold J_Z(X) + \frac 1 {4}\bold J_Z
\bold J_{Z^*}(X),
\endgathered
\tag 1.8
$$
where $X;X^*;Y \in \bold v$ and $Z;Z^*;Z_1;Z_2;Z_3 \in \bold z$
are elements of the Lie algebra.

 By introducing $H(X,X^*,Z,Z^*):=\<\bold J_Z(X),\bold J_{Z^*}(X^*)\>$,
for the Ricci curvature
we have
$$
\gathered
Ricc(X,X^*)=-\frac 1 {2} \sum_{\alpha =1}^lH(X,X^*,e_\alpha ,e_\alpha )
\quad ;\quad Ricc(X,Z)=0; \\
Ricc(Z,Z^*)=\frac 1 {4} \sum_{i=1}^k
H(E_i,E_i,Z,Z^*).
\endgathered
\tag 1.9
$$

By
(1.2), the volume form $\varepsilon$,
the metric tensor
$$
g_{ij}=g\big(\pa_i,
\pa_j\big)\quad ,\quad g_{i\alpha}=g\big( \pa_i,
\pa_\alpha\big)\quad ,\quad g_{\alpha\beta}=g\big(\pa_\alpha
,\pa_\beta\big)
$$
\noindent and the inverse
metric tensor $g^{ij},g^{i\alpha},g^{\alpha\beta}$ have the following
explicit forms:
$$
\varepsilon =dx^1\wedge dx^2\wedge \dots \wedge dx^k\wedge dz^1\wedge
\dots \wedge dz^l,
\tag 1.10
$$
$$
\gather
g_{ij}=\delta_{ij} + \frac 1 {4} \<[X,\pa_i],[X, 
\pa_j]\>=\delta_{ij} + \frac 1 {4}\sum_{\alpha =1}^l
\<\bold J_\alpha\big(X\big),\pa_i\>\<\bold J_\alpha
\big(X\big),\pa_j\>;
\\
g_{i\alpha}=-\frac 1 {2} \<\bold J_\alpha\big(X\big),
\pa_i\> \quad ;\quad g_{\alpha\beta}=\delta_{\alpha \beta},
\tag 1.11
\endgather
$$
$$
\gathered
g^{ij}=\delta_{ij}\quad ; \quad g^{i\alpha}=\frac 1 {2} \<\pa_i,
\bold J_\alpha\big(X\big)\>;
\\
g^{\alpha \beta}=\delta_{\alpha \beta} + \frac 1 {4} \<\bold J_\alpha
\big(X\big),\bold J_\beta\big(X\big)\>.
\endgathered
\tag 1.12
$$

Finally we describe the Dirac- and Laplace-operators acting on forms.

First we 
introduce some minor changes in the notations. So far, we 
used the coordinates $x^i$ and $z^\alpha$ on the Lie algebra.
Next we keep the coordinates $x^i$,while we change
$z^\alpha$ to $x^{k+\alpha}$. 
These new coordinates are denoted by $x^A$,
where $A=1,\dots,k+l$. In the following
 we use also the Einstein summation convention.

Consider a p-form
$$
 \omega=\omega_{(A_1\dots A_p)}dx^{A_1}\wedge\dots\wedge dx^{A_p},
\tag 1.13
$$
where $(A_1\dots A_2)$ implies $A_1<\dots<A_p$.
The $(n-p)$-form $*\omega$ (where $n=k+l$) is defined as usual by
$$
 *\omega={\omega^*}_{(B_1\dots 
B_{n-p})}dx^{B_1}\wedge\dots\wedge dx^{B_{n-p}} ,
\tag 1.14
$$
where
$$
{\omega^*}_{B_1\dots B_{n-p}}=\omega^{(C_1\dots C_p)}
\varepsilon_{(C_1\dots C_p)B_1\dots B_{n-p}}.
\tag 1.15
$$

Then we have:
$$
\gathered
\big(\delta\omega\big)_{A_1\dots A_{p-1}}=
(-1)^{np+n+1}\big(*d*\omega\big)_{A_1
\dots A_{p-1}}=
\\
-g^{AB}\delta_{AA_1 \dots A_{p-1}}^{(B_1 \dots B_p)}\nabla_B
\omega_{(B_1 \dots B_p)},
\endgathered
\tag 1.16
$$
where $\delta_{A_1\dots A_p}^{B_1\dots B_p}=
det\big(\delta_{A_i}^{B_j}\big)$ is
the multivariable Kronecker symbol. 
\smallskip

The Non-Relativistic Dirac- and the
 Laplace-operator (acting on forms) are defined by
$$
d-\delta \text \quad {and} \quad \Delta =(d-\delta)^2=
-(d\delta + \delta d)
\tag 1.17
$$
respectively. Then,
for the Laplacian we have the following well known formula:
$$
\gathered
(\Delta \omega)_{A_1\dots A_p}=
g^{AB}\nabla_A\nabla_B\omega_{A_1\dots A_p}
-\sum_{s=1}^p\omega_{A_1\dots A_{s-1}BA_{s+1}\dots A_p}Ricc^B_{A_s} \\
-\frac 1 {2} \sum_{s=1}^p\sum_{r=1}^p\omega_{A_1
\dots A_{r-1}BA_{r+1}\dots
A_{s-1}AA_{s+1}\dots A_p}{R^{AB}}_{A_rA_s}.
\endgathered
\tag 1.18
$$
\smallskip

We conclude this paragraph by describing the isometries 
on 2-step nilpotent
metric Lie groups.
\proclaim {Proposition 1.1 (\cite {K} \cite {E} \cite {W}) }
The 2-step nilpotent metric Lie groups $\big(N,g\big)$ and $\big(N',g'
\big)$ are isometric iff there exist orthogonal transformations 
$A : \bold v \to
\bold {v'}$ and $C:\bold z \to \bold {z'}$ such that
$$
A \bold J_Z A^{-1} = \bold {J'}_{C( Z)}
\tag 1.19
$$
holds for any $Z \in \bold z$.
\endproclaim

(First Kaplan established this Statement on H-type groups. Then Wilson
generalized it to nilpotent Lie groups. 
Eberlein noticed, that Kaplan's proof
can be applied to 2-step nilpotent case with no altering.
See this Statement also in \cite {GW3} as Proposition 1.4.)

As an example we mention, that the metric groups $\bold H_3^{(a,b)}$ 
(introduced in 
(0.3)) are isometric to the groups $\bold H_3^{\prime(a,b)}$,
where $\bold {J'}_Z$ is
defined by
$$
\gathered
\bold {J'}_Z \big(X_1,\dots ,X_a,X_{a+1},\dots X_{a+b}\big) \\
=\big(X_1Z,\dots ,X_aZ,-X_{a+1}Z,\dots ,-X_{a+b}Z\big).
\endgathered
\tag 1.20
$$
In
this case $C=id$ and
$$
A\big(X_1,\dots ,X_{a+b}\big)=\big(X_1,\dots ,X_a,\overline X_{a+1},\dots
,\overline X_{a+b}\big),
\tag 1.21
$$
where $\overline X$ means conjugation, establish an 
isometry between the two spaces.

 \head
\titlebf 2. \S\quad The general 
$\sigma$- and $\sigma^{(a+b)}$-deformations
\endhead

Let $\{\bold n=\bold v\oplus\bold z,\<,\> ,\bold J_Z\}$ be
a 2-step nilpotent Lie algebra and let 
$\sigma$ be an involutive orthogonal 
transformation on the X-space, commuting with all of the endomorphisms 
$\bold J_Z \in \bold {J_z}$. Then the $\sigma$-deformation of this 
 Lie algebra is
defined by the system
$\{\bold n_\sigma=\bold v\oplus\bold z,\< ,\>,\bold J_Z^\sigma\}$, where
$$
\bold J_Z^\sigma=\sigma\bold J_Z.
\tag 2.1
$$
We call $\bold n$ and $\bold n_\sigma$ also 
{\it $\sigma$-equivalent} Lie algebras.
 Notice that only the endomorphisms
$\bold J_Z$ (i.e. the Lie brackets) are different in the new system,
while the metric structure remains
the old one.

 All
the possible $\sigma$-deformations of $\bold n$ can be easily described
as follows.

Let $\bold v=\bold v_-\oplus \bold v_+$ 
be an orthogonal direct sum such that
both $\bold v_-$ and $\bold v_+$ 
are invariant with respect to the action
of endomorphisms $\bold J_Z\in \bold {J_z}$. 
The corresponding decomposition of
a vector X is described in the form $X=X_- + X_+$. Then 
$\sigma\big(X\big)=-X_- + X_+$ is an appropriate transformation for the 
$\sigma$-deformation
of $\bold n$. If all the invariant subspaces in the
irreducible decomposition of the X-space by the action of the Endomorphism
Space have greater dimension then 2, then
the above construction describes all the $\si$-deformations on the
Endomorphism Space.

The $\sigma^{(a+b)}$-deformations are special $\sigma$-deformations,
which (cf. (0.3) and (1.20)) can be considered as a
more straightforward generalization
of the concept developed for introducing the Lie algebras 
$\bold h^{(a,b)}_3$.

 Notice that in this special case the left products 
$L_Z : \bold H \to \bold H$ (as
well as the right products $R_Z : \bold H \to \bold H$) by imaginary
quaternions $Z$ define an irreducible representation of the Lie algebra
$\bold z=\bR^3=so(3)$ on the 4-dimensional space $\bold H=\bR^4$.

 For generalizing this concept, 
consider an $l$-dimensional linear subspace
$\bJ_{\bold z} \subseteq so(n)$ of orthogonal
endomorphisms. Such a subspace can be 
defined by an appropriate linear map 
$\bJ : \bold z=\bR^l \to End(\bR^n)$. The space $\bJ_{\bold z}$ 
 may or may not be closed with respect to the commutator bracket of 
endomorphisms. 
In the special case of $h^{(a,b)}_3$ it corresponds to the
space $L_{\bR^3}\subset  so(4)$ (or to the space 
$R_{\bR^3}\subset so(4)$), i. 
e.
$\bJ=L$ or $\bJ=R$, where $L$ and $R$ are 
described above.

For a pair $(a,b)$ of natural numbers, the space $\bold v$
is introduced
by the $(a+b)$-times Cartesian product $\bold v=\bR^n \times \dots
\times \bR^n$. The endomorphisms $\bold J^{(a,b)} :\bold z=\bR^l \to
End(\bold v)$ are defined by
$$
\gathered
\bold J^{(a,b)}_Z \big(X_1, \dots ,X_{a+b}\big) \\
=\big(\bold J_Z(X_1), \dots
,\bold J_Z(X_a),-\bold J_Z(X_{a+1}), \dots ,-\bold J_Z(X_{a+b})\big).
\endgathered
\tag 2.2
$$
The corresponding Lie algebra defined by these 
endomorphisms is denoted by
$\bold n^{(a,b)}_J$.

If the space $J_{\bold z}$ is a Lie algebra, then so is the
space $\bold J^{(a+b,0)}_{\bold z}$ and these two 
Lie algebras are isomorphic.
However, the space $\bold J^{(a,b)}_{\bold z}$ with 
$ab\not =0$ can not form
a Lie algebra if $J_{\bold z}$ is a non-Abelian Lie algebra.

In fact,
the involution
$$
\sigma \big(X_1,\dots ,X_{a+b}\big)=
\big(X_1,\dots ,X_a,-X_{a+1},\dots ,-X_{a+b}
\big)
\tag 2.3
$$
 (establishing the deformation $\sigma\bold J^{(a,b)}_\bold z=
\bold J^{(a+b,0)}_\bold z$) commutes
with the elements of both Lie algebras, therefore we have:
$$
\gathered
[\bold J^{(a,b)}_\bold z,\bold J^{(a,b)}_\bold z]=
[\bold J^{(a+b,0)}_\bold z,
\bold J^{(a+b,0)}_\bold z]\subseteq\bold J_{\bold z}^{(a+b,0)}, \\
[\bold J^{(a,b)}_{\bold z},\bold J^{(a+b,0)}_{\bold z}]\subseteq 
\bold J^{(a,b)}_{\bold z}.
\endgathered
\tag 2.4
$$
Since $\bold J_{\bold z}^{(a,b)}\cap
\bold J_{\bold z}^{(a+b,0)}=0$, the space
$\bold J_{\bold z}^{(a,b)}$ can not form a Lie algebra.

In this latter case the space $\bold s^{(a,b)}=
\bold J^{(a,b)}_{\bold z}\oplus
\bold J^{(a+b,0)}_{\bold z}$ forms a 
Lie algebra isomorphic to $\bJ_{\bold z}
\oplus \bJ_{\bold z}$. We will show that 
$\bold s^{(a,b)}$ and $\bold s^{(a',b')}$
with $a+b=a'+b'$ and $ab\not =a'b'$ 
are two inequivalent (non-conjugate) 
representations of the Lie algebra $\bJ_{\bold z}\oplus \bJ_{\bold z}$. This
is the heart of the following Non-Isometry Proof.
\proclaim {Theorem 2.1}
{\bf (A)}
The Lie algebras $\bold n^{(a,b)}_J$ with the same $(a+b)$
and $J$ are $\sigma$-equivalent. Furthermore, the metric Lie group
$N^{(a,b)}_J$ is
isometric to $N^{(b,a)}_J$.

{\bf (B)}
If $\bJ_{\bold z}$ is a non-Abelian Lie algebra 
(or more generally, it contains
a non-Abelian Lie subalgebra), then
$N^{(a,b)}_J$ is locally non-isometric to
$N^{(a',b')}_J$ unless $(a,b)=(a',b')$ up to an order.
\endproclaim
\demo{Proof}
{\bf (A)}The above spaces are clearly 
$\sigma$-equivalent. Furthermore (cf.
 Proposition 1.1), $A=id$ and $C=-id$ establish an isometry between 
$N^{(a,b)}_J$ and
$N^{(b,a)}_J$.

{\bf (B)} First we suppose that $\bJ_{\bold z}$ is a Lie algebra.

If $ab\not =0$, the spaces $\bold J_{\bold z}^{(a,b)}$ and 
$\bold J_{\bold z}^{(a+b,0)}$ can not be conjugate (see (2.4)).
Therefore the metric
Lie groups
$N^{(a,b)}_J$ and
$N^{(a+b,0)}_J$ are non-isometric. In fact, 
if they were isometric, the corresponding
spaces of endomorphisms would be conjugate, by Proposition 1.1. 

To complete the proof in the general case, consider the derived algebra
$\bJ_\bold {z'}=[\bJ_\bold z,
\bJ_\bold z]$ along with the corresponding algebra
 $s^{\prime(a,b)}=\bold J^{(a,b)}_\bold {z'} 
\oplus
\bold J^{(a+b,0)}_\bold {z'}$.  
 Then we have:
$$
s^{\prime(a,b)}=\bold J^{(a,0)}_{\bold z'}\oplus
\bold J^{(0,b)}_{\bold z'} \simeq
\bJ_{\bold z'}\oplus \bJ_{\bold z'}.
\tag 2.5
$$
I. e. $s^{\prime(a,b)}$ can be considered as 
a representation of the Lie algebra
$\bJ_{z'}\oplus \bJ_{z'}$. Let $\chi$ be the character of the corresponding
group representation,
 belonging to the Lie algebra representation $\bJ_{\bold z'}$.
Then the character of the group representation belonging to 
the Lie algebra
representation $s^{\prime(a,b)}$ is 
$\chi^{(a,b)}(g_1,g_2)=a\chi (g_1)+b\chi (g_2)$.
Therefore the representations 
$s^{\prime(a,b)}$ and $s^{\prime(a',b')}$ are inequivalent
 unless $(a,b)=(a',b')$ up to an
order.
 This proves the 
rest part of the Theorem. In fact, by Proposition 1.1, if $N^{(a,b)}_J$
and $N^{(a',b')}_J$ are isometric, 
the spaces $\bold J^{(a,b)}_{\bold z}$ and
$\bold J^{(a',b')}_{\bold z}$ (and consequently also
 the Lie algebras $s^{\prime(a,b)}$ and
$s^{\prime(a',b')}$) are conjugate (see (2.4)).

If $\bJ_{\bold z}$ is not a Lie algebra but it contains a non-Abelian Lie 
subalgebra, apply the above considerations to a maximal non-Abelian Lie
subalgebra in order to have a proof also in these more general cases.
\qed
\enddemo

A different proof of this Statement can be given 
by means of Proposition 5.4,
 where we prove that the Riemannian
 curvatures on $\si$-equivalent spaces are just isotonal
in general.
\smallskip
{\bf Remark 2.2} The proof of the above Theorem breaks down if
$\bJ_{\bold z}$ is just a linear space and it is not closed with respect
to the Lie bracket (commutator) of endomorphisms. 
The following example shows
that the above Theorem can not be generalized to
 such general cases straightforwardly.

In this example $\bold z=\bR^3$, where $\bR^3$ is considered as
 the space of imaginary quaternions.
The endomorphisms $\bJ_Z$ acting on the quaternionic space
$\bold H=\bR^4$ are defined by
$$
\bJ_Z\big(Q\big)=(L+R)_Z\big(Q\big)=ZQ + QZ .
\tag 2.6
$$
(It should be mentioned, that the Lie algebra 
$\big\{\bold n=\bR^4 \oplus
\bR^3,\<,\>,\bold J=R + L\big\}$ was introduced first in \cite {GW3}.)

Since
$$
[(L+R)_{Z_1},(L+R)_{Z_2}]=(L-R)_{Z_1Z_2-Z_2Z_1},
\tag 2.7
$$
the space $(L+R)_\bold z$ is not closed with respect to the Lie bracket
of endomorphisms. Furthermore, for the conjugation $K(Q)=\overline {Q}$
we have:
$$
K(L+R)_ZK^{-1}(Q)=\overline {Z\overline {Q} + 
\overline {Q}Z}=-(L+R)_Z(Q).
\tag 2.8
$$
This means, that for any pair $(a,b)$ of natural numbers the orthogonal
transformations $C=id$ and
$$
A\big(X_1,\dots ,X_{a+b}\big)=\big(X_1,\dots ,X_a,\overline {X}_{a+1},
\dots ,\overline {X}_{a+b}\big)
\tag 2.9
$$
establish an isometry between $N^{(a,b)}_J$ and $N^{(a+b,0)}_J$.
\smallskip
{\bf Remark 2.3}  $\si$-equivalent spaces have inequivalent Riemannian 
curvature in general. For instance, in \cite {Sz3}
 we prove 
(cf. formulas 
(1.17)-(1.20))
that the spaces 
$\bH^{(a,b)}_3$ are distinguished by
the eigenvalues of
 the curvature operator ${R_{ab}}^{cd}$.
 (See more about this problem in \S 5.)

 However, by (1.8), the $\si^{(a+b)}$-equivalent spaces have 
identical Ricci curvature at the origin, providing new examples for
non-isomorphic Lie groups with locally non-isometric invariant metrics
yet having equivalent Ricci curvature.

The reader can find such examples
in \cite {Ka}. The main difference between these two type of examples is
that in Karidi's examples the spaces have different degrees of
the balls' volume growth, while in our examples the spaces have the same
Riemannian density function on polar coordinate neighborhood around the
origin (this latter statement will be proven in a subsequent paper).

 \head
\titlebf 3. \S\quad Technicalities on 
Ball$\times$Torus- and on Ball-type
domains
\endhead

The spectrum coincidences will be considered on the following two types
of manifolds:

(1) Let $\Gamma$ be a full lattice on the Z-space spanned by a basis
$\{e_1,\dots ,e_l\}$. 
For an $l$-tuple $\alpha =(\alpha_1,
\dots ,\alpha_l)$ of integers the corresponding lattice point is 
$Z_\alpha =\alpha_1e_1+
\dots 
+\alpha_le_l$. Since $\Gamma$ is a discrete subgroup, 
one can consider the
factor manifold $\Gamma \backslash N$ with the factor metric. This
factor manifold is a principal fiber bundle with the base space 
$\bold v$
and with the fiber $T_X$ at a point $X\in \bold v$. Each fiber $T_X$ is
 naturally identified with 
the torus $T=\Gamma \backslash 
\bold z$. The projection 
$\pi : \Gamma \backslash N \to \bold v$ defined by
$\pi : T_X \to X$ projects the inner product from 
the horizontal subspace
($=$ orthogonal complement 
 to the
fibers) to the euclidean inner product $\<,\>$ on the X-space.

Consider
also a euclidean ball $B_\delta$ of radius $\delta$ 
around the origin of the
X-space and restrict the fiber bundle onto $B_\delta$. Then this fiber
 bundle
$(B_\delta ,T)$ has the boundary $(S_\delta ,T)$, 
which is also a principal
fiber bundle over the sphere $S_\delta$.

When we consider the spectrum problems on these manifolds, we call these
cases $Ball\times Torus$- and $Sphere\times Torus$-cases.

Only
these manifolds with the Laplacian spectrum on functions
 were investigated in the literature
 so far.

(2)In this paper we consider also such domains around the 
origin which are homeomorphic
to a $(k+l)$-dimensional ball, and their smooth boundaries can be 
described as 
level sets by equations of the form
 $f(|X|,Z)=0$. In this case the boundary is homeomorphic to
a sphere $S^{k+l-1}$ such that 
the boundary points form
a euclidean sphere of radius $\delta (Z)$,
for any fixed Z. 
 I. e. the boundary can be
described by the equation $|X|^2-\delta^2 (Z)=0$. We
call these cases {\it Ball-cases} resp. {\it Sphere-cases}.

 In the following we need also the explicit form of
 the normal vector $\bmu$ at the boundary of these domains.
From the equation
$$
\nabla f=grad f=\sum_i\bold X_i(f)\bold X_i+
\sum_\al \bold Z_\al (f)\bold Z_\al
$$
we get (by using the special function 
$f(|X|,Z)=|X|^2-\delta^2 (Z)$) that this unit normal vector at a
point $(X,Z)$ is
$$
\bmu =(4|X|^2+\frac 1 {4} |\bold J_{\nabla\delta^2}(X)|^2+
|\nabla\delta^2|^2)^
{-\frac 1 {2}}(2X-\frac 1 {2}\bold J_{\nabla\delta^2}(X)-
\nabla\delta^2),
\tag 3.1
$$
where $\bmu$ is considered as an element of the Lie algebra. 
Notice that in
Case $(1)$ the $\bmu$ has the simple form $\bmu =X/|X|$.

 It can be written in the following regular vector form
  (cf. (1.2)): 
$$
\bold \bmu=A\big(2|X|
E_0 - \sum_{\al =1}^l(\pa_\al \delta^2)
(\frac 1 {2}|X|E_\al +\sum_{\beta=1}^l(1+\frac 1 {4}\<\bold J_\al (X),
\bold J_\beta (X)\>)e_\beta )\big ),
\tag 3.2
$$
where $\{e_1,\dots ,e_l\}$ is an orthonormal basis of $\bold z$ and
$E_\al =\bold J_\al (E_0);E_0=X/|X|$ furthermore
$$
A=(4|X|^2+\frac 1 {4} |\bold J_{\nabla 
\delta^2}(X)|^2 + |\nabla\delta^2|^2)^{-\frac 1 {2}}.
\tag 3.3
$$

In the literature 
this is the first paper which 
involves also 
the manifolds of type $(2)$ into the
spectral investigations. The considerations 
require completely different techniques
on these two types of manifolds.  
 
 Since $(\bold J^\si_Z)^2=
\bold J^2_Z$, the endomorphisms 
$\bold J_Z$ and $\bold J^\si_Z$ are isospectral
(conjugate). Therefore the following
 Statement immediately follows from the Gordon-Wilson
Theorem \cite {GW3}. However, we describe a proof 
here which is a straightforward
adoption of the proof given on the 
special spaces $\bold H^{(a,b)}_3$ in \cite {Sz3}. The main difference
between this proof and the proof given in \cite{GW3} is that our proof
does not use the total geodesic assumption on the tori of the
considered torus bundle. This idea was used also in \cite{GS}.

\proclaim {Theorem 3.1} The metrics on
 $\si$-equivalent spaces are Dirichlet and Neumann isospectral
on any of the manifolds described in the Ball$\times$Torus-cases. 
\endproclaim

\demo {Proof}
First consider the Weierstrass decomposition
$$
L^2(\Gamma\backslash N)=\oplus W^{(\alpha)}
\tag 3.4
$$
of the $L^2$ function space, where $W^{(\alpha)}$ consists of functions 
of the form
$$
\phi^{(\alpha)}(X,Z)=\varphi (X)e^{2\pi \bold i\<Z_\al ,Z\>},
\tag 3.5
$$
for any fixed lattice point $Z_\alpha$.

By (1.6), any of the spaces $W^{\al}$ is invariant under 
the action of the
Laplacian. More precisely we have:
$$ 
\Delta \phi^{(\al)}=(\square_{(\al)}\varphi )e^{2\pi 
\bold i\<Z_\al ,Z\>},
\tag 3.6
$$
where
$$
\square_{(\al)}=\Delta_X + 2\pi \bold i D_{(\al )}\bullet -4\pi^2
\big(|Z_\al|^2 + \frac 1 {4} |\bold J_{Z_\al}(X)|^2\big).
\tag 3.7
$$

Next we construct intertwining operator, $\tau_*$, which intertwines the
corresponding
Laplacians.

Since $\bold J^\si_{Z_\al}$ and $\bold J_{Z_\al}$ are 
isospectral (conjugate),
there exists an orthogonal transformation $\tau^{(\al)}$ 
on the X-space such that
$$
\bold J^\si_{Z_\al}=\tau^{(\al)}\circ \bold J_{Z_\al}\circ 
\tau^{(\al)-1}.
\tag 3.8
$$
 (In \cite {Sz3}, this transformation $\tau^{(\al)}$
 is individually constructed,
using conjugations and products by $Z_\al$.) Therefore we have
$$
D^\si_\al\bullet =\tau^{(\al)}_*D_\al\bullet\tau^{(\al)*}\quad ; \quad
\square^{\si}_{(\al )} = \tau^{(\al)}_*\square_{(\al)}\tau^{(\al)*}.
\tag 3.9
$$
Thus the operator $\tau_* =\oplus \tau^{(\al)}_*$ 
intertwines the Laplacians
$\Delta^\si$ and $\Delta$ on the whole $L^2$ function spaces.

It is plain that $\tau_*$ does not change 
the Dirichlet boundary conditions.
By (3.1),
the normal vector
is $\bmu =X/|X|$ 
 at a boundary point $(X,Z)$. 
 Since the
orthogonal transformations leave the differentiations from the radial
directions invariant, also the Neumann boundary 
condition  is invariant by the
action of $\tau_*$.
This proves the Statement completely.
\qed
\enddemo

 The above applied method of Fourier transform breaks down in the
Ball-case, since the functions of the form (3.5) do not satisfy
neither the Dirichlet nor the Neumann boundary conditions. 
For constructing an appropriate intertwining operator in the Ball-cases,
 we have to 
develop a completely different technique.

{\bf Remark 3.2} Any 2-step nilpotent Lie group $N$ 
can be extended into a
solvable group $SN$ defined on the half space $\bold n\times \bR_+$ with
multiplication given by
$$
(X,Z,t)(X',Z',t')=(X+t^{\frac 1 {2}}X',Z+tZ'+\frac 1 {2}t^{\frac 1 {2}}
[X,X'],tt').
\tag 3.10
$$
This group has the Lie algebra $\bold s=\bold n\oplus \bold t$ 
with the Lie 
bracket satisfying
$$
[\pa_t,X]=\frac 1 {2}X\quad ;\quad [\pa_t,Z]=Z\quad;\quad 
[\bold n,\bold n]_{/SN}
=[\bold n,\bold n]_{/N},
\tag 3.11
$$
for any $X \in \bold v$ and $Z \in \bold z$.
 Consider on this Lie algebra
the natural inner product $\<,\>$ such that we keep the inner product on
$\bold n$ furthermore $\pa_t\perp \bold n$ and $|\pa_t|=1$. 
The left invariant
extension of this inner product is denoted by $g$.

The explicit computations of the several geometric objects along
with other technicalities on 
these solvable extensions are established in \cite {Sz4}.
For instance, we prove there that two appropriate solvable extensions
 are 
isometric iff
the corresponding nilpotent groups are isometric (cf. Proposition 1.2
in \cite{Sz4}; this Statement is established also in \cite{GW1}). 
The Laplacian acting on
functions has the following explicit form (cf. (1.12) in \cite{Sz4}):
$$
\gathered
\Delta=t\Delta_X+t^{\frac 1 {2}}\Delta_Z+\frac 1 {4}t\sum_{\al ;
\beta =1}^l
\<\bold J_\al (X),\bold J_\beta (X)\>\pa^2_{\al \beta} \\
+t\sum_{\al =1}^l\pa_\al D_\al\bullet+t^2\pa^2_t+(1-\frac k 
{2}-l)t\pa_t.
\endgathered
\tag 3.12
$$
(This formula is established also in \cite{GS}.)

The above spectral investigations can be extended onto 
these solvable groups
on the following domains.

{\it The Solvable Ball $\times$ Torus Cases:} A group $SN$ can be 
considered as a
principal fiber bundle (vector bundle) over the $(X,t)$-space, fibered by
the Z-space. Again, $\Gamma$ is a full lattice on the Z-space,
 while $D_{R(t)}$
is a domain on the $(X,t)$-space diffeomorphic to a $(k+1)$-dimensional
ball whose smooth boundary can be described by an equation of the form
$|X|=R(t)$. Then consider the 
torus bundle 
$(\Gamma /\bold z,D_{R(t)})$ 
over $D_{R(t)}$.
The normal vector $\bold \bmu$ at a boundary point
 $(X,Z,t)$ is of the form
$\bold \bmu=A(t)X+B(t)\pa_t$, where $A(t)$ and $B(t)$ 
are determined by $R(t)$.

{\it The Solvable Ball Case:} In this case we consider a domain 
$D$ diffeomorphic
to a $(k+l+1)$-dimensional ball whose smooth boundary (diffeomorphic to
$S^{k+l}$) can be described  as a level surface in the form
 $|X|=R(Z,t)$.
The normal vector $\bold \bmu$ at a boundary point 
$(X,Z,t)$ has the form
$$
\gathered
\bmu = A_0(|X|,Z,t)E_0
+C(Z,t)\pa_t + \\
\sum_{i=1}^l(A_i(|X|,Z,t)E_i+B_i(|X|,Z,t)e_i),
\endgathered
\tag 3.13
$$
where the coefficient functions 
$A_\al B_i$ and $C$ are determined by $R(Z,t)$ (cf. (3.24)-(3.26) in
\cite{Sz4}).

Then, by repeating the very same arguments used in the nilpotent cases
 we get
(the only modification needed is that
the function $\varphi$, in (3.5), should have the
form $\varphi (X,t)$) 
\proclaim {Theorem 3.3} The metrics on
 $\si$-equivalent solvable spaces are Dirichlet and Neumann
isospectral on any of the domains described in the Solvable 
Ball$\times$Torus cases.
\endproclaim
{\bf Remark 3.4} Finally, we check that how these ideas
are (or, are not) working in the case of forms.

In the Ball$\times$Torus-cases,
the Fourier decomposition on the $L^2$ space of r-forms is
$$
L^2_r(\Gamma /N)=\oplus W^\al_r, 
\tag 3.14
$$
where $W^\al_r$ is spanned by the r-forms of the form
$$
\Omega^{(\al)}=\varphi (X)e^{2\pi \bold i\<Z_\al,Z\>}\omega.
\tag 3.15
$$
Here $\varphi (X)$ is an $L^2$-function and $\omega$ is an invariant(!)
r-form.

The invariant r-forms can be expressed by the invariant 1-forms
$$
\gathered
\theta_i(.)=g(\bold X_i,.)=dx^i(.), \\
\vartheta_\beta (.)=g(\bold Z_\beta,
.)= dz_\beta (.)-\frac 1 {2}\sum_{i=1}^l
\<\bold J_\beta (X),\pa_i\>dx_i(.).
\endgathered
\tag 3.16
$$
For computing the right sides of these formulas we used 
(1.1) and (1.11). 
From these formulas we immediately get that the forms $\theta_i$ are
closed, while 
$$
d\vartheta_\beta =-\frac 1 {2}\bold J_\beta^{ji}\theta_j\wedge \theta_i.
\tag 3.17
$$

Actually this non-closedness of the forms $\vartheta_\beta$ 
"kills" the method
of the Fourier transform if we try to prove the
isospectral property of $\sigma$-equivalent spaces on the spectra
of forms of $degree\geq 1$. In fact, in case of functions 
the method of Fourier transform
reduces the spectral investigation to a case when only a
 single endomorphism
$\bold J_\al$ is involved. The above equation shows that such kind of
reduction is not possible in case of higher order forms. (Check 
 \S 5. to see that the same problem emerges also with respect to the
$\delta$-operator.)

 \head 
\titlebf 4.\S\quad  
Constructing Isospectral Balls on Spaces with Anticommutators
\endhead 

The final goal of this section is to find such/those 
particular $\si$-deform\-ations
providing isospectral metrics also on the Ball- 
and Sphere-type manifolds.
It turns out that this goal can be achieved on spaces having particular
Endomorphism Spaces described as follows.

 If we perform $\si$-deformations
 only on a particular subspace of the Endomorphism Space
$\bold {J_z}$ and we keep the endomorphisms belonging to the
 orthogonal complement unchanged, such a deformation is said to be 
{\it partial $\si$-deformation} or {\it $P\si$-deformation}. 
More precisely, if the
deformation is performed on the endomorphisms belonging to a
subspace $\bold s\subset \bold z$, the deformation is called 
$\si_{\bold s}$-deformation and if $\si$ has the form $\si =(id,-id)$
with respect to a decomposition 
$\bold v=\bold v^{(a)}\oplus\bold v^{(b)}$,
it is called  {\it $\si_{\bold s}^{(a+b)}$-deformation}.
 
If all the invariant
subspaces of the irreducible decomposition of the 
X-space with respect to the
action of $\bold {J_z}$ have dimension greater than 2, then a 
$\si_{\bold s}$-deformation is a 
$\si^{(a+b)}_{\bold s}$-deformation. For the sake of the simplicity, 
we suppose that the Endomorphism Spaces satisfy 
 this {\it non-degenerating property}.
By this assumption
 only trivial cases are excluded (these cases are really trivial when
The Endomorphism Space contains an anticommutator).

 These deformations define a much
 wider class then the $\si$-deformations.
 However, a $P\si$-deformation
may be equivalent to a $\si$-deformation.  
In the following we prove
that a 
$\si$-deformation is equivalent to a 
$\si_{\bold A}$-deformation 
if the endomorphism space $\bold{J_z}$ contains a
so called anticommutator $A$.
(Here 
$\bold A$ denotes the 1-dimensional subspace
spanned by $A$.) It turns out that 
these particular $\si\simeq\si_{\bold A}$-deformations are the desired
ones and they provide
isospectral metrics also on the Ball- and 
Sphere-type manifolds and not only on
the Ball$\times$Torus- and the Sphere$\times$Torus-type manifolds.
However, it is still an open question, whether
 they provide all the solutions for the
problem.

We establish all these Theorems in 3 Steps.

\medskip\noindent
{\bf Endomorphism Spaces with Anticommutators (alias $ESW_A$)}
\medskip

A non-degenerated endomorphism $A=\bold J_Z$ is said to be 
an {\it anticommutator} 
in $\bold {J_z}$  
if
$A\circ B=-B\circ A$ holds for all $B\in\bold J_{Z^\perp}$, 
i. e., if the $A$  
anticommutes with all of the orthogonal endomorphisms. An 
anticommutator satisfying $A^2=-id$ is said to be a 
{\it unit anticommutator}.

 Any anticommutator can be written in the
form $A=S\circ A_0$, where $A_0$ is a unit endomorphism anticommuting
with the elements of $\bold A^\perp$ and the symmetric scaling
endomorphism $S$ (which is one of the square root of $-A^2$) 
is commuting with all elements of the $ESW_A$.

An orthogonal transformations $\sigma$,
 used for defining $\si$-deformations on a non-degenerated Endomorphism
Space, is always
 of the form $\si =(id,-id)$ with respect to an appropriate
 decomposition
$\bv =\bv^{(a)}\oplus\bv^{(b)}$.
In these cases, the $\si$-deformation of an endomorphism $\bJ_Z$ is
denoted by $\bJ^{(a,b)}_Z$ and the original endomorphism can
be denoted also by
$\bJ_Z^{(a+b,0)}$. 

Since the equation
$ 
\bold J_Z\bold J_{Z^*}+\bold J_{Z^*}\bold J_Z=-2\<Z,Z^*\>id
$ holds 
on any Heisenberg type group 
(cf. (1.4) in \cite {CDKR}),
all the endomorphisms are anticommutators in the
Endomorphism Spaces belonging to these groups.

In this paper we consider such
Endomorphism Spaces where the
anticommuting elements form
 at least a 1-dimensional subspace.
 The complete
description of the Endomorphism Spaces having non-trivial anticommuting 
elements is given 
in \cite{Sz4} (cf. Canonical Matrix Representation Theorem 2.1
of the $ESW_A$'s). Here we introduce only the Quaternionic $ESW_A$'s,
which show the wide range of such spaces.

Let $A$ be a $k\times k$ diagonal matrix 
such that it has the same imaginary 
quaternionic number $\bold a$ in the main diagonal.
 Consider also a 
linear space $\bold J^\perp _A$ of $k\times k$ symmetric matrices 
such that the entries are imaginary quaternions perpendicular to 
$\bold a$.
All these matrices define a skew-symmetric real
 matrix-map on $\bold H^k=\bR^{4k}$.
 (In this matrix transformation we multiply by 
the entries of the matrix 
always on the same side,
i. e. we can choose between left or right representations 
of these matrices.)
Then the endomorphism $\bold J_A$ is an anticommutator 
in the endomorphism
 space $\bold J_A^\perp \oplus \bold A$. Also the
Cartesian products of such spaces are endomorphism spaces with
non-trivial
anticommutators.

Surprisingly enough, if an endomorphism space 
has an anticommutator $A$, a
$\si$-deformation is equivalent to the $\si_{\bold A}$-deformation.
 More precisely we have
\proclaim {Reduction Theorem 4.1}  
If
$\bJ_{\bold s}$ is an odd-dimensional subspace
 of the anticommutators in a non-degenerated $ESW_A$, then 
a $\si =\si^{(a+b)}$-deformation is equivalent to the partial
 $\si_{\bold s}=\si_{\bold s}^{(a+b)}$-deformation,
i. e. the corresponding metrics 
are isometric. 

If the space $\bold s$ is even dimensional, then a
 $\si_{\bold s}=\si_{\bold s}^{(a+b)}$-deformation
  provides a space which is
 isometric to the original space.
\endproclaim
\demo{Proof}
Let $\bold {J_z}$ be a linear space of endomorphisms 
with
 a non-degenerated anticommutator
$\bJ_A$. 
By an appropriate rescaling of the maximal invariant
spaces of $\bJ_A$,
 one can find an endomorphism $A_0$ (which may not be in $\bold{J_z}$) 
such that $A^2_0=-id$.   
 Since all the endomorphisms $B\in \bJ_{A^\perp}$ 
leave the maximal eigensubspaces
of $A$ invariant, 
 this $A_0$ is commuting with $A$ and it is anticommuting with the 
endomorphisms from the orthogonal space $J_{A^\perp}$.

Let $\{J_{A_1},\dots ,J_{A_r}\}$ be an orthonormal basis on the subspace
$J_{\bold s}$. The orthogonal transformation $A_{i0}^{(b)}$ on the space
$\bold v=\bold v^{(a)}\oplus\bold v^{(b)}$ is 
defined by $id\times
 A_{i0/\bv^{(b)}}$. 
 Then the orthogonal transformation $A_0^{(b)}=A_{10}^{(b)}
\circ\dots\circ A_{r0}^{(b)}$ and $C=id$ establish the isometries
(cf. Proposition 1.1)
between the corresponding spaces.

In fact, if $r=2m+1$, we have
$$
(A_0^{(b)})^{-1}\bold J_Z^{(a+b,0)}A_0^{(b)}=\bold J_Z^{(a,b)},
\forall Z\in\bold s^\perp
\quad ;\quad
(A_0^{(b)})^{-1}\bold J_{A_i}^{(a,b)}A_0^{(b)}=\bold J_{A_i}^{(a,b)}.
$$
If $r=2m$, then we get
$$
(A_0^{(b)})^{-1}\bold J_Z^{(a+b,0)}A_0^{(b)}=\bold J_Z^{(a+b,0)},
\forall Z\in\bold s^\perp
\quad ;\quad
(A_0^{(b)})^{-1}\bold J_{A_i}^{(a,b)}A_0^{(b)}=\bold J_{A_i}^{(a+b,0)},
$$
which prove the statement completely.
\qed
\enddemo

In the following we consider an endomorphism space $J_{\bold z}$ 
 with an anticommutator $J_A$. For the sake of the simplicity 
we suppose that $J_A^2=-id$, i. e.
it is a unit anticommutator. 
 All of the
theorems proven in the following can be established also for
general anticommutators. In the end we describe all the modifications
what should be done for having the proofs in these general cases.

The great advantage on the $ESW_A$'s is that
 a $\si$-deformation can be expressed as a
$\si_{\bold A}$-deformation and an intertwining operator
for the corresponding Isospectrality Theorem can be established by
means of the endomorphism $A$.
  
\medskip\noindent
{\bf Harmonic Analysis based on a Unit Anticommutator}
\medskip 

 The main tool in this paper is a sort of
 Harmonic Analysis established by means of a unit anticommutator.
We define a certain decomposition of the $L^2$ function space 
with the help of this unit endomorphism and then we use
this decomposition for constructing appropriate intertwining operators. 

For a fixed X-vector
$Q$, we define the complex valued function
$$
\bold \Theta_Q(X)=\<Q+\bold i\bold J_A(Q),X\>.
\tag 4.1
$$
Then the complex valued function space is spanned by 
the functions of the form
$$
\Phi_{(Q_i,p_i,Q^*_i,p^*_i,\varphi )} (X,Z)=\big(\Pi_{i=1}^r
\bold \Theta_{Q_i}^{p_i}(X)
\overline {\bold \Theta}_{Q^*_i}^
{p^*_i}(X)\big)\varphi (Z).
\tag 4.2
$$

Notice, that the functions of the pure form
$$
\Phi_{(Q_i,p_i)}(X)=\Pi_{i=1}^r\bold 
\Theta_{Q_i}^{p_i}(X)\quad ;\quad \overline{ 
\Phi}_{(Q_i,p_i)}(X)=\Pi_{i=1}^r\overline {\bold \Theta}_{Q_i}^{p_i}(X)
\tag 4.3
$$
are harmonic with respect to the euclidean Laplacian $\Delta_X$ 
on the $X$-space. However, the polynomials of the mixed form are not
harmonic, since 
$$ 
\Delta_X\Theta_Q(X)\overline Q_{Q^*}(X)=2\<Q,Q^*\>+2\bold i
\<\bold J_A^{(a,b)}(Q),Q^*\>.
\tag 4.4
$$ 
 In the following we find 
more specific functions which still span the $L^2$
function space on the X-space. First we rewrite  
 the
functions (4.2) in the form
$$
\varphi (|X|,Z)\Phi_{(Q_i,p_i,Q^*_i,p^*_i)}(X_u),
\tag 4.5
$$
where $X_u=X/|X|$ is the unit
radial X-vector.

Notice that
  the  $q^{th}$ order
polynomials
$$
\gathered
\Phi_{(Q_u,q)} (X)=\<Q_u+\bold i \bold J_A(Q_u),X\>^q, \\ 
\overline {\Phi}_{(Q_u,q)}(X)=
\<Q_u -\bold i \bold J_A(Q_u),X\>^q
\endgathered
\tag 4.6
$$
are harmonic on the euclidean X-space whose restrictions 
onto the unit sphere $S$ 
(around the origin of the X-space) define $q^{th}$
order eigenfuctions of the Laplacian $\Delta_S$. Let us note that
 the whole space 
$\bold H^{(q)}$ of $q^{th}$ order eigenfunctions is not spanned
by these functions. 
The "missing" functions can be furnished with the help
of the $q^{th}$ order functions of the form 
$\Phi_{(Q_i,p_i,Q^*_i,p_i^*)}
(X_u)$, where
$q=\sum_i p_i+\sum_ip_i^*=s+(q-s)$, as follows

 Such a function is not an eigenfunction
of the Laplacian $\Delta_S$, however, it is
 an eigenfunction
of the differential
operator $D_A\bullet$ with the eigenvalue $(2s-q)\bold i$.  
By orthogonal projection, project all these functions 
onto the function space 
$\bold H^{(q)}$.
The image space is denoted by $\bold H^{(s,q-s)}$.
 
 In \cite {Sz3} we proved that the operators $\Delta_S$ and
$D_A\bullet$ commute, where $D_A\bullet$ means 
directional derivative
with respect to the vector field $X\to \bold J_A(X)$.
(This proof is based on the simple fact that the vector field
$J_A(X)$ tangent to the X-spheres around the origin is an infinitesimal
generator of isometries on these spheres.) 
Therefore the whole space $\bold H^{(q)}$ is invariant under
the action of $D_A\bullet$ and the subspace 
$\bold H^{(s,q-s)}\subset \bold H^{(q)}$
is an eigensubspace of this operator with the eigenvalue 
$(2s-q)\bold i$.
Thus the decomposition 
$\bold H^{(q)}=\oplus_{s=0}^q\bold H^{(s,q-s)}$ is an 
orthogonal direct sum corresponding to 
the common eigensubspace decomposition
of the two commuting differential operators $\Delta_S$ and 
$D_A\bullet$.

We need a more accurate description of the functions which span these
eigensubspaces. For this description we use the 
  kernel functions 
$$
H_{(q)}(Q_u,Q^*_u)=\sum_{j=1}^{N_q}\eta^{(q)}_j(Q_u)
\eta^{(q)}_j(Q^*_u),
$$ 
where $\{\eta^{(q)}_1,\dots ,\eta^{(q)}_{N_q}\}$
 is an orthonormal basis on the subspace $\bold H^{(q)}$.
 In \cite {Be} it
is proven (cf. Lemma 6.94) that
for any fixed $Q_u$, the eigen
function $H_{(q)}(Q_u,.)\in \bold H^{(q)}$ is a radial one 
(i.e. it has the form
$C_q\<Q_u,.\>^q+\dots +C_1\<Q_u,.\> +C_0$) with
$H_{(q)}(Q_u,Q_u)=1$.  
(One can find appropriate reference to this Statement also
in \cite {Sz1} (\S 2, Pages 8-11) where this kernel is scrutinized
 on arbitrary harmonic spaces.)
 These kernel functions can be used for projecting 
the functions into
the subspaces $\bold H^{(q)}$. 
 More precisely, for any function $\Psi\in L^2(S)$ we have 
$$
\Psi(Q_u)=\sum_{q=0}^\infty\int H_{(q)}(Q_u,Q^*_u)\Psi (Q^*_u) dQ^*_u=
\sum h_{(q)}(\Psi )(Q_u),
\tag 4.7 
$$
   which is called 
the {\it spherical
decomposition of $\Psi$ by the spherical harmonics}. The operators
$
h_{(q)} :L^2(S)\to\bold H^{(q)}
$
project the $L^2$ function space to the corresponding eigensubspace of 
the Laplacian. 

One of the most important properties of these operators is
  the commutativity
with the differential operators $D_\al\bullet$. This commutativity
immediately follows from the fact that the $\bold J_\al (X)$ is an 
infinitesimal generator of one parametric families of 
isometries on the euclidean 
sphere $S$, furthermore, the projections $h_{(q)}$ are invariant
with respect to these isometries.

For an X-vector $Q$, we define the complex vector 
$$
\bold Q=Q+\bold i\bold J_A(Q).
\tag 4.8
$$ 
Then,
by substituting $Q_u=1/2(\bold Q_u+
\overline {\bold Q}_u)$ into the above
expression of the radial kernel $H_{(q)}(Q_u,.)$, we get
\proclaim {Proposition 4.2} The $r^{th}$ order polynomial space $\bold
P^{(r)}$ is the direct
sum 
of the subspaces $\bold P^{(p,r-p)}$ spanned 
by the polynomials of the form
$\Theta_{Q}^p\overline {\Theta}_{Q}^{r-p}$, where $Q\in\bold v$
and $0\leq p\leq r$. The space 
$\bold P^{(p,r-p)}$ consists of the $r^{th}$
order eigen polynomials 
of the differential
operator $D_A\bullet$
with eigenvalue $(r-p)\bold i$.

The projection 
$h_{(r)}$
 establishes a one to one and onto map between the polynomial spaces
 $\bold P_S^{(p,r-p)}$ and $\bold H^{(p,r-p)}$.
 We denote this restricted
map by 
$h_{/(r)}: \bold P_S^{(r)}\to\bold H^{(r)}$.
 The direct sum 
$T=\oplus_rh_{/(r)}$ of these maps
 defines an invertible operator on the whole function space 
$L_S^2$. This operator $T$ is commuting with the differential operators
$D_\al\bullet$.
\endproclaim 

Above
$\bold P_S^{(p,r-p)}$ denotes the space of 
 the corresponding restricted functions on the sphere $S$. 
This space is spanned by the functions of the form 
$\Theta_{Q_u}^p
\overline {\Theta}_{Q_u}^{r-p}$, 
where $Q_u$ is a unit vector. Let us mentioned
that the projected functions 
$h_{/(r)}(
\Theta_{Q_u}^p
\overline {\Theta}_{Q_u}^{r-p})$ 
can be written in the following form
$$
h_{/(r)}(\Theta_Q^p\overline{\Theta}_Q^{r-p})=\sum_{s=0}^{min(p,r-p)}
A_s\Theta_Q^{p-s}\overline{\Theta}_Q^{r-p-s} ,
$$
where $A_0=1$ furthermore the higher order coefficients can be computed
by the recursive formula
$$
2|Q|^2(p-s)(r-p-s)A_s+(s+1)(dim(\bold v)+r-2)A_{s+1}=0.
$$

This statement follows from the fact that the homogeneous polynomial
$\sum_s A_s\<X,X\>^s
\Theta_Q^{p-s}\overline{\Theta}_Q^{r-p-s}$
is harmonic with respect to the Laplacian $\Delta_X$ and therefore
its restriction onto $S$ is an eigenfunction of the Laplacian
$\Delta_S$. Actually this statement proves the stronger Proposition
asserting that the projected functions onto the lower order subspaces
have the form
$$
h_{(q)}(\Theta_Q^p\overline{\Theta}_Q^{r-p})=\sum_{s=0}^{min(p,q-p)}
A_{qs}\Theta_Q^{p-s}\overline{\Theta}_Q^{q-p-s} ,
$$
where $q\leq r$ and they are trivial $(=0)$ on the higher order 
subspaces.
 
By (4.7), the radial eigenfunctions span the eigenspaces
 $\bold H^{(r)}$,
therefore the above formulas completely describe the maps $h_{/(r)}$.
However, the above concept can be easily generalized to the case when
$h_{/(r)}(\Pi_i
\Theta_{Q_i}^{p_i}
\overline {\Theta}_{Q_i}^{r_i-p_i})$ should be figured out, where
$\sum r_i=r$. 
In this case we look for the harmonic homogeneous polynomial
of the form
$$
\Pi_i\sum_{s_i} A_{s_i}\<X,X\>^{s_i}
\Theta_i^{p_i-s_i}\overline{\Theta}_i^{r_i-p_i-s_i}.
\tag 4.9
$$

 To compute clean
recursion formulas in this general cases is tiresome and fortunately
we do not need them in this paper. However, the recursive
computation of the coefficients $A_{s_i}$ is clear. Arrange this latter
function according to the $s$-lines (the 0-line corresponds to the
original function) and then let the Laplacian
act on this function. First the coefficients $A_{1_i}$ can be
computed by considering the action of the Laplacian on the 0- and the
1-line. E. t. c., the computation should be continued till it is
finished.

This Computational Technique will be used later in this section.

These formulas together with $h_{/(r)}D_\al\bullet =
 D_\al\bullet h_{/(r)}$ completely prove the second part of
the Proposition.

In the following we represent the functions $\psi\in\bold H^{(r)}$
in the form $\psi=h_{/(r)}(\psi^*)$, where $\psi^*\in\bold P^{(r)}_S$.
We use this tool to simplify the considerations.

On 
the whole ambient space, $\bold n$, the function
space is spanned by the functions of the form
$$
F(X,Z)=\varphi (|X|,Z)h_{/(r)}(\Theta^p_{Q_u}
\overline {\Theta}_{Q_u}^{r-p})(X_u),
$$
where $X_u=X/|X|$ is a unit vector. The decomposition with respect to
these 
functions is called {\it spherical decomposition}. 
If we cancel $h_{/(r)}$
in (4.9), the corresponding decomposition 
with respect to these functions
is called {\it polynomial decomposition}. The operator $T$ transforms
a polynomial decomposition to a spherical decomposition.
\medskip\noindent
{\bf Constructing the Intertwining Operator}
\medskip

If $\bv =\bv^{(a)}\oplus\bv^{(b)}$ is the decomposition corresponding
to a $\si_A$-deformation on a space with non-degenerated $ESW_A$,
 then the $\si_A$-deformed group is denoted by $N_{\bz}^{(a,b)}$.
The original group $N_{\bz}$ can be denoted also by $N_{\bz}^{(a+b,0)}$.
This general notation is consistent with the particular notation
$\bold H^{(a,b)}_3$, however, the $\si^{(a+b)}$-deformation is
represented now as $\si_A^{(a,b)}$-deformation. If 
$N=N^{(a,b)}_{\bold z}$ and 
$N'=N^{(a',b')}_{\bold
z}$, where $a+b=a'+b'$, are $\si_A$-deformations of the very same group,
then $N$ and $N'$ are clearly
 $\si_{\bold A}$-equivalent spaces. Consider these groups, the 
Intertwining Operator $\kappa$ between the corresponding complex 
function spaces is defined by 
$$
\gathered
\kappa :L^2(N_{\bold z}^{(a,b)},\bold C)\to 
L^2(N_{\bold z}^{(a',b')},\bold C),\\
\kappa :
F(X,Z)\to
F'(X,Z)=\varphi (|X|,Z)h'_{/(r)}
(\Theta{\prime p}_{Q_u}\overline{\Theta '}^{r-p}_{Q_u})
(X_u),
\endgathered
\tag 4.10
$$
where $F$ is defined above by means of $\bJ_A^{(a,b)}$
 and the functions $\Theta '$ are
 defined on $N'$ by means of $\bold
J_A^{(a',b')}$.

Notice that this map maps the function space $\bold H^{(r)}$ onto
$\bold H^{\prime(r)}=\bold H^{(r)}$,
 however, it is defined by means of the map
$$
\kappa^*:\bold P^{(r)}\to\bold P^{\prime(r)}\quad ,\quad
\kappa^*:
\Theta^p_Q\overline{\Theta}_Q^{r-p}\to
\Theta^{\prime p}_Q\overline{\Theta}_Q^{\prime r-p}
$$
and of the projection $h_{(r)}$. Unlike $\kappa$, 
the map $\kappa^*$ can be easily handled.
If $\{E_1,\dots ,E_K\}$ is a basis on the X-space then
$\kappa^*(\Pi\Theta_{E_{i_r}}\Pi\overline{\Theta}_{E_{j_r}})=
\Pi\Theta'_{E_{i_r}}\Pi\overline{\Theta}'_{E_{j_r}}$
and in general, for arbitrary vectors $Q_m$, we get
$$
\kappa^*(\Pi\Theta_{Q_{i_r}}\Pi\overline{\Theta}_{Q_{j_r}})=
\Pi\Theta'_{Q_{i_r}}\Pi\overline{\Theta}'_{Q_{j_r}}.
$$

The operators $\kappa$ and $\kappa^*$ are connected by the equation
$\kappa=T'\circ\kappa^*\circ T^{-1}$.
 In the following we show 
the intertwining property of the map $\kappa$. 

Consider the explicit expression
(1.6) of the Laplacian $\Delta$.
Since $\kappa :\bold H^{(q)}\to \bold H^{\prime(q)}=\bold H^{(q)}$, 
the terms $\Delta_X + \Delta_Z$ and $\Delta '_X +\Delta '_Z$ are 
 intertwined by this map. 
In fact, the intertwining of $\Delta_Z=\Delta^{\prime}_Z$ is obvious
and the intertwining of $\Delta_X=\Delta^{\prime}_X$ follows from the
decomposition $\Delta_X=\pa^2_{|X|}+(dim\bold v-1)|X|^{-1}\pa_{|X|}
+|X|^{-2}\Delta_S$.

Choose an orthonormal basis $\{e_0,e_1,\dots ,e_{l-1}\}$ on the Z-space
such that $e_0=A$. The Greek characters are used for the indices
$\{0,1,\dots ,l-1\}$ and the Latin characters are used for the indices
$\{1,\dots ,l-1\}$. Then
$D_c \bullet =D'_c \bullet$,
furthermore
$$
D_0\bullet \bold \Theta_Q(X)=\bold i\bold \Theta_Q(X)\quad ,\quad 
D_0\bullet \overline {\Theta}_Q(X)=-\bold i\overline {\Theta}_Q(X),
\tag 4.11
$$
and
$$
D_c\bullet
\Theta_{Q}=\<\bold Q,\bold J_c^{(a,b)}(X)\>=
-\overline{\Theta}_{\bold J_c(Q)}\quad ,\quad
D_c\bullet\overline{\Theta}_{Q}=-
\Theta_{\bold J_c(Q)}.
\tag 4.12
$$
(Let us mention that the switching of the conjugation in (4.12)
is due to the equation $\bold J_0\bold J_c=-\bold J_c \bold J_0$).
Therefore 
$\kappa^*D_\al\bullet (\psi^*)
=D'_\al\bullet
\kappa^*(\psi^*)$
for any $\psi^*\in\bold P^{(r)}_S$. 
Since $\kappa =T'\circ\kappa^*\circ T^{-1}$, we also have
$\kappa D_\al\bullet (\psi)
=D'_\al\bullet
\kappa(\psi)$ and thus
the terms $\sum \pa_\al D_\al \bullet$ and $\sum \pa_\al D'_\al \bullet$
 in (1.6)
are intertwined by the map $\kappa$.

Only the term $(1/4)\sum \<\bold J_\al (X),\bold J_\beta (X)\>
\pa_{\al \beta}$ should be considered yet.

First notice that on the Heisenberg-type groups, particularly on
$\bold H^{(a,b)}_3$, this
operator is nothing but $(1/4)|X|^2\Delta_Z$. Therefore the 
Intertwining Property of $\kappa$ is perfectly
established in these special groups and (later we see that) also
on their solvable extensions. Particularly, the Intertwining Property
is established for the Striking Examples and only the intertwining
of the boundary conditions should be scrutinized to proving the
Isospectrality Theorem.

We turn back to the investigation of the considered
 operator on general spaces.
 Since $\bold J_0\circ 
\bold J_c$ is a skew symmetric endomorphism,
 $\<\bold J_0(X),\bold J_c(X)\>=0$. Therefore this 
 operator is the same on the 
$\sigma_{\bold A}$-equivalent spaces $N$ and $N'$.
Yet, for the sake of completeness we should prove that the $\kappa$ maps
a function of the form 
$$
\<\bold J^{(a+b,0)}_c(X),
\bold J^{(a+b,0)}_d(X)\>h_{(r)}(\Theta^p\overline
{\Theta}^{r-p}):=\bold J_{cd}(X)h_{(r)}(\Theta^p\overline{\Theta}^{r-p})
\tag 4.13
$$ 
to a function of the very same form on $N'$. 

This problem can be 
simplified by the Spherical Decomposition Theorem (4.7). By this 
Theorem, the radial spherical harmonics span all the eigensubspaces
$\bold H^{(r)}$, therefore, it is enough to consider the real functions
of the form 
$K_{(r)}(X)=J_{cd}(X)h_{/(r)}\<Q,X\>^r$. 
In this case the function
$h_{/(r)}\<Q,X\>^r$ is nothing but a constant multiple of the function
$H_{(r)}(Q,X)$.

To figure out $\kappa (K_{(r)})$, first the spherical decomposition
of $K_{(r)}$ should be computed. By the "Technique" establishing the
coefficients $A_{s_i}$ in (4.9), the component spherical harmonics
are built up by the functions of the form
$$
\gathered
\<Q,X\>^p\, ,\,
J_{cd}(X)\<Q,X\>^q\, ,\\
 \<(J_cJ_d+J_dJ_c)(Q),X\>^s
\<Q,X\>^{v}=\<Q_{cd},X\>^s  
\<Q,X\>^{v},
\endgathered
\tag 4.14
$$  
 such that the combinational coefficients depends only on
the constants $r,p,s,v,TrJ_c\circ
J_d$ and on $\<\bold J_c(Q),\bold J_d(Q)\>$.

The function $\bold J_{cd}(X)$ should be further scrutinized. It
can be written in the form
$$
\gathered
\bold J_{cd}(X)=
\sum_{i=1}^K \frac 1 {4} \{ (
\< X,\bold Q_{ci}\>\< X,\overline {\bold Q}_{di}\> +\< X,\overline
{\bold Q}_{ci}\>\< X,\bold Q_{di}\> ) + \\ (
\< X,\bold Q_{ci}\>\< X,\bold Q_{di}\> +\< X,\overline
{\bold Q}_{ci}\>\< X,\overline{\bold Q}_{di}\> )
\} = \bold J^{(1)}_{cd}(X)+\bold J^{(2)}_{cd}(X),
\endgathered
\tag 4.15
$$
where $E_1,\dots ,E_K$ is an orthonormal basis on the X-space 
($K=k(a+b)$) and
$\bold Q_{ei}=\bold J_e^{(a,b)}(\bold E_i)$. 
The proof of this formula immediately follows from 
$$
\bold J_{cd}(X)=\sum_i\<\bold J_c(X),E_i\>\<\bold J_d(X),
E_i\>\quad ,\quad E_i
=\frac 1 {2} (\bold E_i+\overline{\bold E}_i).
\tag 4.16
$$

Notice that the second function in (4.15)
is a spherical harmonic function on the unit sphere $S$, 
having a vanishing integral
on $S$. This statement follows from formulas 
$\<\bold Q_1,\bold Q_2\>=
\<\overline{\bold Q}_1,\overline{\bold Q}_2\>=0$.
 By the very same reason,
also the other two integrals vanish in 
the spherical expansion of this function.
I. e., the function $\bold J^{(2)}_{cd}(X)$
 vanishes and the simplification $\bold J_{cd}(X)=\bold
J^{(1)}_{cd}(X)$ can be introduced in formula (4.15). 

By the substitution $Q=(1/2)(\bold Q+\overline{\bold Q})$
we get that the component-spherical-harmonics in the 
decomposition of the function $K_{(r)}$ 
are the linear combinations of the functions of the form
$$
\Theta^{p-s}_{Q}\overline{\Theta}^{r-p-s}_{Q}\, ,\,
 \bold J_{cd}
\Theta^{p-s}_{Q}\overline{\Theta}^{r-p-s}_{Q}\, ,\,
(\Theta_{Q_{cd}}+\overline{\Theta}_{Q_{cd}})
\Theta^{p-v}_{Q}\overline{\Theta}^{r-p-v}_{Q}
\tag 4.17
$$
 such that the combinational coefficients depends only on
the constants $r,p,s,v,TrJ_c\circ
J_d$ and on $\<\bold J_c(Q),\bold J_d(Q)\>$.

The considered problem can be settled by representing $J_{cd}$ in the
form (4.15) in formula (4.17).
In fact, the coefficients discussed above
 are the same on both spaces (since they
do not depend on the unit anticommutator $\bJ_0^{(a,b)}$). By
 the Technique establishing the coefficients $A_{s_i}$ in formula
(4.9) we get that
 the pre-images 
(with respect  to the map $T$) 
of the considered component-spherical-harmonics
are
combinations of
 appropriate 
 functions
(gotten from (4.17)) such that the coefficients do not depend on
$\bJ_0^{(a,b)}$. This completely proves
that the $\kappa$ maps a function $K_{(r)}$ (or, even a more 
complicated function, described (4.13)) to an appropriate function
desired in this problem.

To complete the proof of the Isospectrality Theorem, we have to prove
yet that the intertwining map $\kappa$ does not change the boundary 
conditions.

Since the functions of the form $\varphi (|X|,Z)$ are mapped to the
very same function by the $\kappa$, furthermore, in any 
of the considered cases
the Dirichlet boundary condition imposes restrictions only on these
functions, the $\kappa$ obviously preserves the 
Dirichlet boundary conditions.
Also the Neumann boundary condition is preserved.

In fact, the coefficient
functions in
  (3.2) are the same  
on $\si$-equivalent spaces. Furthermore, 
by (4.11) and (4.12), the $\kappa$ intertwines the directional 
derivatives with respect to the radial vectors $X$ as well as with
respect to the vectors $\bold J_{\al}(X)$. Therefore 
$\nu (\phi)=0$ if and only $\nu '\phi '=0$.

By summing up, we have
\proclaim {Main Theorem 4.3} 
The operator $\kappa=T'\circ\kappa^*\circ T^{-1}$
intertwines the corresponding Laplacians and boundary conditions. 
Therefore the metrics on
 the considered $\si_{\bold A}$-equivalent spaces 
$N^{(a,b)}_{\bold z}$
 and 
$N^{(a',b')}_{\bold z}$ with $a+b=a'+b'$ 
are both Dirichlet
and Neumann isospectral on any of the regions described in the Ball-
or Ball$\times$Torus-cases.
\endproclaim

According to Theorem 2.1, these constructions provide a wide range of 
isospectral spaces with different local geometries. Among these examples
we find the H-type group $\bold H^{(a,b)}_3$. 
In this
case
the endomorphism algebra $J_{\bR^3}$ forms a Lie algebra. All the 
endomorphism spaces containing this Lie algebra such that, say $\bJ_1$
is an anticommutator in the whole space 
(the wide range of these possibilities
are described at the beginning of this paragraph!) 
also provide such examples.

It should be mentioned that
 these constructions provide the first examples for pairs
of locally non-isometric yet Dirichlet and Neumann isospectral manifolds
constructed on the most simple simply connected domains 
 diffeomorphic to euclidean balls.
 
\noindent
{\bf Remark 4.4} All the above statements can be easily 
extended to the more
general cases, when $\bJ_A$ is just a non-degenerated anticommutator.
 In this case, first,
 introduce the unit anticommutator $A_0$ by an appropriate 
rescaling of $A$.
This $A_0$ may not be in the endomorphism space, however, 
it is commuting
with $A$ and it is anticommuting with the elements of 
$\bold J_{A^\perp}$.
Then establish $\kappa$ by using $A_0$. Notice that
 the function $\bold J_{00}(X)$ 
 is not 
equal to the function $|X|^2$ in general. However, 
formulas (4.14) and (4.15)
remain true in the same form with respect to this function, which proves
the statement completely. 

Actually the Theorem can be stated in a more 
general form:
\proclaim {Main Theorem $4.3^{\prime}$}
 Let $\bold{J_z}=\bold A\oplus\bold A^\perp$ and 
 $\bold{J_z^\prime}=\bold A^\prime\oplus\bold A^{\prime\perp}$ 
be endomorphism
 spaces with anticommutators such that the endomorphisms are acting
on the same space $\bold v$ furthermore $\bold A^\perp=\bold A^{\prime
\perp}$ and the anticommutators $A$ and $A^\prime$ are isospectral
(conjugate). Then the Ball-type domains are
both Dirichlet and Neumann
isospectral with respect to the metrics $g$ and $g'$
 of the spaces $N_J$ and $N_{J^\prime}$.
\endproclaim

This more general Theorem is proven in \cite{Sz4}.
By this more general Theorem,
one can construct isospectral spaces by 
deforming the anticommutator $A$ isospectrally,
 while keeping the orthogonal complement $\bold A^\perp$
unchanged. The full exploration of this idea is given in \cite{Sz4}.

This exploration includes the complete description of the $ESW_A$'s in a
Matrix Representation Theorem as well as the complete description of
all possible non-trivial isospectral deformations of an anticommutator.
It turns out that the continuous deformations are always trivial,
providing isometric spaces. However, beside the 
$\si_{\bold A}$-deformations we point out a cornucopia of 
deformations providing spectrally inequivalent Endomorphism Spaces.
They represent a completely new type of examples ever constructed
in this field. For instance, these metrics are isospectral only
on the Ball- and Sphere-type domains and they are not isospectral
on the Ball$\times$Torus-type domains or on their boundaries. 

\noindent
{\bf Remark 4.5}  Using (3.12) and 
(3.13), the above Isospectrality Theorem can be trivially generalized 
to the solvable spaces $SN^{(a,b)}_{\bold z}$ and the Theorem
remains true in the very same form. 

In both cases the Theorem can be
 established also on the boundaries of these domains. 
Extremely interesting
examples are provided by the geodesic spheres with the same
radius on the harmonic spaces
$S\bold H^{(a,b)}_3$,
 where $a+b$ is the same dimension number. It turns out that these
spheres are isospectral and one of them is homogeneous while the others
are inhomogeneous. We consider these Examples in \S 6.

 \head
\titlebf 5.\S\quad Spaces with isotonal curvatures. 
Some results on invariant
forms  
\endhead

We study the analogous isospectrality properties 
of the Laplace- and the non-relativistic Dirac-spectra on forms in
 a separate paper. Here we describe some
results only about the spectra of the Laplacian acting 
on invariant differential forms. In the end
we consider also the spectrum of the curvature operator acting as a 
symmetric endomorphism on the 2-forms. We prove that 
$\si$-equivalent spaces have isotonal
curvature. Both type of these spectra have an important impact on the
corresponding
spectra of the general forms. 

We start by describing the action of $d$ and $\delta$ 
operators on invariant
differential forms.

The invariant differential forms can be expressed by the 1-forms
$\theta_i$ and $\vartheta_\al$ introduced in (3.16),
since all these 
differential forms can be expressed as a linear combination
of the forms:
$$
\theta_{\{i_1,\dots ,i_{p_1},\al_1,\dots ,\al_q\}} = 
\bold \theta_{i_1}\wedge\dots\wedge\theta_{i_{p_1}}\wedge
\vartheta_{\al_1}\wedge\dots\wedge\vartheta_{\al_q},
\tag 5.1
$$
where $i_1<\dots <i_{p_1}$
and $\al_1<\dots <\al_q$.

The formulas (3.17) completely describe the action of the d-operator on 
invariant differential forms. Therefore we should deal only with the
$\delta$ operator.
From (1.16) we straightforwardly get
$$
\gathered
\delta \theta_i=0\quad ,\quad \delta\vartheta_\al =0\quad ,\quad
\delta (\theta_i\wedge\vartheta_\al )=0, \\
T_{ij}=\delta (\theta_{i}\wedge\theta_{j})
=-\vartheta_{[E_{i},E_{i}]}=-
\sum_{\al=0}^l\<\bold J_\al (E_{i}).E_{j}\>\vartheta_\al, \\
\delta (\theta_{i_1}\wedge\dots\wedge\theta_{i_p}\wedge\vartheta_{\al_1}
\wedge\dots\wedge\vartheta_{\al_q})= \\
\sum_{r\<s,}C
\dots\wedge\theta_{i_{r-1}}\wedge
\theta_{i_{r+1}}\wedge\dots\wedge\theta_{i_{s-1}}\wedge\theta_{i_{s+1}}
\wedge\dots\wedge\vartheta_{\al_q}\wedge T_{i_ri_s},
\endgathered
\tag 5.2
$$
where C is a combinatorial constant depending on $p$ and $q$.

The Laplacian acting on invariant differential forms can be expressed
by these formulas.

 In the following we describe some results about the spectra of the
 Laplacian acting on invariant differential forms. Such a description 
of the
spectra on the invariant forms 
is an important step also to the 
description of the spectra on the general 
forms.
 In fact, the general forms are spanned by the forms 
$F\theta_{\{i_1,\dots ,\al_q\}}$, and
the action of the Laplacian on such a form is
$$
\Delta F\theta_{(i_1,\dots ,\al_q)}
=(\Delta F) \theta +2\nabla_{grad F}\Theta
+F\Delta\theta .
$$
Once the spectrum problem is solved on functions and on the invariant 
differential forms (i. e. on the ends of the above formula), the general
problem requires the investigation of the middle operator 
$\nabla_{grad F}
\theta$. 

 In this paper we determine the spectrum of the Laplacian only
 on invariant 1- and 2-forms.

Since the $\delta$-operator is trivially acting on 1-forms,
$\Delta =-\delta d$ in this case. By (3.17) and (5.2) we get
\proclaim {Proposition 5.1}
The eigenvectors of the Laplacian acting on the invariant 1-forms are
divided into two classes.

(A) The X-forms $\theta_Q$ with eigenvalue 0.

(B) The Z-forms $\vartheta_{Z_\al}$ such that $Z_\al$ is an eigenvector
of the quadratic form $T(Z,Z^*)=Tr \bold J_Z\circ\bold J_{Z^*}$. The
eigenvalues are the negatives of the eigenvalues of this quadratic form.

Therefore, 
  the Laplacians acting 
on invariant 1-forms  
are isospectral 
on $\si$-equivalent spaces.
\endproclaim

It should be mentioned that this spectrum is 
the same also on the members
of the continuous family of two step nilpotent Lie groups constructed
by Gordon and Wilson in \cite {GW3} (see Corollary 4.2 of this article).
In view of Schueth's Theorem (namely, these spaces are not isospetral
on general 1-forms), the isospetrality on functions and on invariant
1-forms does not imply the 
isospectrality on the general 1-forms. (Though
Schueth's Theorem concerns the boundary of these manifolds, it is
extremely probable that the statement is valid on the manifolds as well.

Next we determine the spectrum of 
the Laplacian acting on the invariant 2-forms.

For the sake of simplicity we suppose 
that the inner product on the Z-space
is the natural inner product 
$T(Z_1,Z_2)=Tr \bold J_{Z_1}\circ\bold J_{Z_2}$.
If we use a different inner product, 
then, like in the previous statement,
the eigenvalues of the quadratic form T
emerge in the following spectrum. 

The space $\wedge^2$ of the 
invariant 2-forms splits into the following 3
invariant subspaces:
$$
\theta_{\bold v}
\wedge\theta_{\bold v}\quad , \quad \vartheta_{\bold z}\wedge
\vartheta_{\bold z}\quad , 
\quad \theta_{\bold v}\wedge\vartheta_{\bold z}.
\tag 5.3
$$
By (3.23) and (5.2), 
the Laplacian is acting 
on these invariant subspaces 
as follows.
$$
\gathered
\Delta \theta_i\wedge\theta_j=-\sum_\al\bold J_{\al ij}\bold J_\al^{pq}
\theta_p\wedge\theta_q\quad ,\quad 
\Delta \vartheta_\al\wedge\vartheta_
\beta=-2\vartheta_\al\wedge\vartheta_\beta
\\
\Delta \theta_i\wedge\vartheta_\al=\theta_i\wedge\vartheta_\al +
\sum_{\beta,j}\<\bold J_\beta (E_i),\bold J_\al (E_j)\> \theta_j\wedge
\vartheta_\beta= \\
\theta_i\wedge\vartheta_\al +\sum_{\beta ,j}(\<\bold J_{\al}(E_i),\bold
J_\beta (E_j)\> +
\<[\bold J_\al ,\bold J_\beta ](E_i),E_j\>)\theta_j\wedge\vartheta_\beta
\endgathered
\tag 5.4
$$

One can interpret these formulas as follows.

(A) In case of the invariant subspace 
$\theta_{\bold v}\wedge\theta_{\bold v}$,
a 2-form $u^{ij}\theta_i\wedge\theta_j$ can be considered as a skew 
endomorphism. Then 
the Laplacian is the $-1$ times of the projection onto
the subspace $\bold J_{\bold z}$ in this case.
 Therefore the Laplacian has the eigenvalue $-1$ 
with multiplicities $l$ and the eigenvalue 0 with multiplicities equal
to the codimension of $\bold J_{\bold z}$ in the whole space of skew 
endomorphisms. 

The second case should not be interpreted.

(B) With respect to the invariant space 
$\theta_{\bold v}\wedge\vartheta_{\bold z}$
notice that the restricted Laplacian onto this space is the same on
$\si$-equivalent spaces. 
Therefore the $\si$-equivalent spaces are isospectral
on the invariant 2-forms.

\proclaim{Proposition 5.2}
The Laplacians acting on the invariant 2-forms
 are isospectral  
on $\si$-equivalent spaces.
\endproclaim

Finally we consider the curvature operators on these spaces. By (1.9)
we obviously get
\proclaim{Proposition 5.3}
Two $\si$-equivalent spaces 
have identical Ricci curvature at the origin.
\endproclaim

 Now we consider the spectrum of the curvature operator ${R_{ij}}^{pq}$
 acting as a symmetric endomorphism on the space of 2-forms. This
spectrum is called the $\wedge^2-spectrum$ of the curvature operator. By
(1.18), this spectrum has a great impact 
on the Laplacian spectra on forms.
 Yet, surprisingly enough, this spectrum distinguishes
 the $\si$-equivalent spaces, since we have
\proclaim {Proposition 5.4} 
The curvatures of the $\si$-equivalent spaces
$N^{(a,b)}_{\bold z}$ and
$N^{(a',b')}_{\bold z}$  
 with $aba'b'\not =0$ are always isotonal (i. e.,
the elements of the spectra
are the same but the multiplicities may be different).

In many cases they are strictly isotonal (not isospectral) 
unless
$(a,b)=(a',b')$ up to an order. For instance, this statement is true on
the spaces $\bold H^{(a,b)}_3$.

A general
 criteria is formulated in the following
proof for this Isotonality Statement: 
If $0\not =ab\not =a'b'\not =0$, the 
curvatures are strictly isotonal 
iff the spectrum on the
mixed boxes changes when it is multiplied by $-1$.

The curvature of the space $N^{(a+b,0)}_{\bold z}$ is subtonal to the 
curvature  
 of the above spaces, i. e. the corresponding tonal spectrum on 
$N^{(a+b,0)}_{\bold z}$ is a subset of the common tonal spectrum of the
other spaces.
\endproclaim
 \demo{Proof} Also in this case we decompose the space $\wedge^2$
 of the 2-forms into invariant
subspaces.
Unlike in case of the Laplace spectrum
on the invariant 2-forms, 
 now the decomposition $\bold v=\bold v^{(a)}\oplus\bold v^{(b)}$,
where $\bold v^{(a)}=\bR^{na};\bold v^{(b)}=\bR^{nb}$ (see (2.2)),
 plays an important role. 
First we decompose the space
 $\wedge^2$ into the following
3 invariant subspaces:
$$
\gathered
\bold D=\theta_{\bold v^{(a)}}\wedge\theta_{\bold v^{(a)}}\oplus
\theta_{\bold v^{(b)}}\wedge\theta_{\bold v^{(b)}}\oplus
\vartheta_{\bold z}\wedge\vartheta_{\bold z} \\ \bold F=
\theta_{\bold v^{(a)}}\wedge\theta_{\bold v^{(b)}}\quad ,\quad \bold G=
\theta_{\bold v}\wedge\vartheta_{\bold z}
\endgathered
\tag 5.5
$$
On the space $\bold G$, 
the action is described by
$$ \theta_i\wedge\vartheta_\al
\to \sum_{j,\beta}-\frac 1 {4}\<\bold J_\beta (E_i),\bold J_\al (E_j)\>
\theta_j\wedge\vartheta_\beta .
\tag 5.6
$$ 
In Proposition 5.2 we proved that this operator is 
 isospectral on any pair of $\si$-equivalent spaces.
 
 On the space $\bold D$, the curvature is acting by 
the components described by the first and
the last lines in (1.8). 
In this case we need a more refined decomposition
into invariant subspaces.

 In the following
$\bold v_r^{(a)}\simeq\bR^n$ 
means the $r^{th}$ component-space of $\bold 
v^{(a)}$.  
The component-space
$\bold v^{(b)}_s$ is defined similarly.  

The forms forming the subspace
$$
\bold {Dg}=\sum_r
\theta_{\bold v_r^{(a)}}\wedge\theta_{\bold v_r^{(a)}}\oplus\sum_s
\theta_{\bold v_s^{(b)}}\wedge\theta_{\bold v_s^{(b)}}\oplus
\vartheta_{\bold z}\wedge\vartheta_{\bold z} 
\tag 5.7
$$
are said to be diagonal forms. 

The space $\bold {Dg^\perp}\subset
\bold D$ of the orthogonal forms in $\bold D$ is spanned by the
blocks 
$\bold v_r^{(a)}\wedge\bold v^{(a)}_{r^*}$ and
$\bold v_s^{(b)}\wedge\bold v^{(b)}_{s^*}$, where $s\not =s^*$
and $r\not =r^*$. Note that on these blocks only the component
$$
\gathered
u^{i_rj_{r^*}}\theta_{i_r}\wedge\theta_{j_{r^*}}\to
\frac 1 {2}\sum_\al u^{i_rj_{r^*}}
{\bold J_{\al i_r}}^{p_r}
{\bold J_{\al i_r^*}}^{p^*_{r^*}}\theta_{p_r}\wedge\theta_{p^*_{r^*}},\\
u^{i_sj_{s^*}}\theta_{i_s}\wedge\theta_{j_{s^*}}\to
\frac 1 {2}\sum_\al u^{i_sj_{s^*}}
{\bold J_{\al i_s}}^{q_s}
{\bold J_{\al i_s^*}}^{q^*_{s^*}}\theta_{q_s}\wedge\theta_{q^*_{s^*}}
\endgathered
\tag 5.8
$$
of the curvature operator (which comes from the second and third term
of the first line of (1.8)) acts non-trivially,
leaving this space
invariant. Furthermore this operator is the same on any of the
considered spaces, therefore the curvature is isospectral on
these invariant spaces. 

We can prove the same isospectrality also on the invariant
space $\bold{Dg}$ of the diagonal forms.

In fact, on a space $N^{(a,b)}_{\bold z}$ consider the map
$\tau :\bold{Dg}\to
\bold {Dg}$ such that 
$\tau =id$ on $\theta_{\bold v^{(a)}}\wedge\theta_{\bold v^{(a)}}$ and
$\tau =-id$ on $\theta_{\bold v^{(b)}}\wedge\theta_{\bold v^{(b)}}$.
Then the $\tau$ intertwines the curvatures of the spaces 
$N^{(a+b,0)}_{\bold z}$ and
$N^{(a,b)}_{\bold z}$ on this invariant space $\bold {Dg}$.

None of the spectrums investigated so far
 distinguishes the considered spaces.
 These spaces are distinguishable by the curvature spectrum
on the invariant space $\bold F$. This space splits into the direct
sum of the "mixed" boxes
$\theta_{\bold v^{(a)}_r\wedge\bold v^{(b)}_s}$.
Each of these boxes is invariant on which only the component
$$
u^{i_rj_s}\theta_{i_r}\wedge\theta_{j_s}\to
\frac 1 {2}\sum_\al u^{i_rj_s}
{\bold J_{\al i_r}}^{p_r}
{\bold J_{\al i_s}}^{q_s}\theta_{p_r}\wedge\theta_{q_s}
\tag 5.9
$$
of the curvature is acting non-trivially. On a pair of such boxes,
the spectrum is the 
same. 
 More precisely,
this spectrum is the negative of the spectrum what we get on the boxes
of the invariant space $\bold{Dg}^\perp$.

If this spectrum does change when it is multiplied by $-1$,
 the curvatures on
the spaces $N^{(a+b,0)}_{\bold z}$ and $N^{(a,b)}_{\bold z}$ can not
be isospectral. 
However, the first is subtonal to the second one by the above
arguments. If $0\not =ab\not =a'b'\not =0$, the spaces $N$ and $N'$
have different number of mixed boxes. Therefore the curvature spectra
of these spaces on the space of mixed boxes are just isotonal.

Let us mentioned that the above mentioned spectrum on a mixed box
of the space $\bold H^{(a,b)}_3$ is
$\{
\lambda =-1/2,m_\lambda =12;
\lambda =3/2,m_\lambda =4\}$, 
 where $m_\lambda$ means the multiplicity
of the corresponding eigenvalue. 
This spectrum changes if we multiply it by
$-1$, therefore these $\si$-equivalent 
spaces have strictly isotonal curvatures (unless 
$(a,b)=(a',b')$ up to an order). 
\qed
\enddemo

 \head 
\titlebf 6.\S\quad
 Constructing the Striking Examples
\endhead

The proofs of the Isospectrality 
Theorems on the boundary manifolds are incorporated
into \cite{Sz4}. We have to point
out that the proof of the isospectrality itself
 is the easy part of these
considerations,
since from the explicit expression of the Laplacian on $\pa D$
 one can get immediately that an appropriate Intertwining Operator can
be defined in the very same way as on the ambient space. Only
the functions $\varphi (|X|,Z)$
(resp. $\varphi (|X|,Z,t)$ in the solvable case)
 should be changed to functions of the
form $\varphi (\delta (Z),Z)=\phi (Z)$ 
(resp. of the form $\varphi (\delta (Z,t),Z,t)=\phi (Z,t)$), 
where $Z$ varies on the boundary. However, the Non-Isometry Proofs are
much more complicated on the boundaries than on the ambient spaces,
 because the boundaries are inhomogeneous
spaces and we lose the technique of the homogeneous 
spaces.

 In the following we describe more details on the special H-type spaces 
$\bold H^{(a,b)}_3$ and on their solvable extensions 
$S\bold H^{(a,b)}_3$,
to discover the most striking examples of our
constructions. Namely, the metrics on the geodesic spheres with the 
same radius are isospectral
on the spaces 
$S\bold H^{(a,b)}_3$ with the same $a+b$, despite the radical 
differences shown by the isometry groups: The isometry groups have
different dimensions on these spaces and, even more surprisingly,
the geodesic spheres are homogeneous on the 2-point homogeneous
space
$S\bold H^{(a+b,0)}_3$, while the metrics on the
 other spheres are locally 
inhomogeneous.

 In the following we consider domains, $D$, such that
the boundary $\pa D$ is described by a function of the form 
$|X|=\delta (|Z|)$ in the nilpotent case and of the form
$|X|=\delta (|Z|,t)$ in the solvable case.
I. e., both the X- and the Z-cross sections with this boundary are spheres
denoted by $S_X(Z)$ and $S_Z(X)$ in the nilpotent case and by
 $S_X(Z,t)$ and $S_Z(X,t)$ in the 
 solvable case.
Let us mention that the geodesic spheres around the origin 
are such hypersurfaces 
on these
spaces. 
More precisely, the geodesic
spheres with the same radius $R$ around the origin
 are the same point-sets on these
groups  and they are described by the same equation which has a
relatively complicated implicit form,
 $\sigma (R,|X|,|Z|)=0$, in the nilpotent case (see 
\cite {CDKR} and \cite {Sz2} for more details) and it has the
more simple explicit form
$$
|X|^2=4((e^R+e^{-R}+2)t-|Z|^2)^{1\over2}-4(t+1)
\tag 6.1
$$
in the solvable case \cite{Sz4}.

The Laplacian on the boundary $\pa D$ can be computed
by standard arguments. The corresponding formulas are established
in \cite{Sz4}. One can check the idea of these computations also
in \cite {Sz3} (cf. formulas (3.1)-(3.3)). In the nilpotent case we
get
$$
\tilde \Delta =\Delta_{S_X(Z)}+\big (1+\frac 1 {4}\delta^2 (|Z|)\big )
\Delta_{S_Z(X)}+ 
\sum_{\al =1}^3\big (\pa_\al 
-Z_{u\al})D_\al\bullet .
\tag 6.2
$$
The last term is nothing but the restriction of the operator
$\sum_\al \pa_\al\bold J_\al\bullet$ onto the boundary $\pa D$.

In the solvable case we have
$$
\gathered
\tilde
\Delta =t\Delta_{S_X}+(t^{\frac 1 {2}}+{1\over4}t|X|^2)\Delta_{S_Z}+ 
 \\
t\sum_{\al =1}^3(\pa_\al -Z_{u\al}) D_\al\bullet+
c^2t^2(\pa_t-\bmu_t)^2 -
c^2(2+\frac k {2})t(\pa_t -\bmu_t),
\endgathered
\tag 6.3
$$
where $\bmu_t=C(|Z|,t)\pa_t$ is introduced in (3.13).

Restrict the functions (4.9) onto the boundary $\pa D$
and define the map $\tilde\kappa$ on these
restricted functions again by the formula (4.10). The same idea is
applicable to the solvable extensions. Then
 by the same proof established on the ambient spaces we get that this
map intertwines the Laplacians on the considered function spaces. 
Therefore we have
\proclaim {Theorem 6.1} If $a+b=a'+b'$, 
 the metrics induced
on the Sphere-type boundary $\pa D$ are isospectral on the spaces 
$\bold H^{(a,b)}_3$ and
$\bold H^{(a',b')}_3$
 as well as on the solvable extensions
$S\bold H^{(a,b)}_3$ and
$S\bold H^{(a',b')}_3$.

 Particularly we get that the metrics on the geodesic spheres 
with the same radius 
of these
spaces
are isospectral.
\endproclaim

The Non-Isometry Proofs can be established by
a simple modification of the Proofs given in \cite {Sz3} 
on the Sphere$\times$Torus-type manifolds. Actually, three
different proofs are given there and each of them can be reconstructed
also for the case of Sphere-type boundaries.

The details are as follows

At a point $(X,Z)\in \pa D$, let $\wt z$ be the tangent space of the
sphere $S_Z(X)$, furthermore consider also the distribution
$\rho^{(a,b)}=\bold J^{(a,b)}_{\bold z}(X)\subset T(S_X(Z))$
on the sphere $S_X(Z)$ of the X-space. In the solvable case we consider
also the 1-dimensional distribution $\wt T$ 
spanned by the unit tangent vector
$\bold t$ of the t-parameter curves on the surface. Let
$K^{(a,b)}$ be the distribution orthogonal to 
$\rho^{(a,b)}\oplus \wt z$ in the
nilpotent case and it should be
 perpendicular also to $\wt T$ in the solvable
case.

By computing the Ricci curvature on $\pa D$ it turns out that the 
$K^{(a,b)}$ is an eigensubspace of the
Ricci operator and the eigenvalues, taken on this subspace, are
 completely different from the other eigenvalues. This statement
concerns the nilpotent case as well as the solvable case.

For proving this Statement, the matrix of the Ricci operator is computed
in a particularly chosen basis in \cite{Sz4}. At a point $(X,Z)$ (resp.
at
$(X,Z,t)$) on the surface the orthonormal
basis $\{\bold i,\bold j,\bold k\}$ on the Z-space is chosen such that
$\bold i=Z_u$. Then $\bold j$ and $\bold k$ are tangent to the surface.
If $\bmu_X$ (considered as a Lie algebra element) is the X-component
of the normal vector and $\bmu_{Xu}$ denotes the normalized vector,
then introduce the vectors (considered as Lie algebra elements)
$$
\wt E_{\bold i}=\bold{J_i}(\bmu_{Xu})\, ,\, 
\wt E_{\bold j}=\bold{J_j}(X_u)\, ,\, 
\wt E_{\bold k}=\bold{J_k}(X_u).
\tag 6.4
$$
Let us mentioned that the vectors $\wt E_{\bold j}$ and 
$\wt E_{\bold k}$ are tangent to the distribution $\rho^{(a,b)}$ while
$\wt E_{\bold i}$ is not tangent to this distribution, expressing the
fact that the subspaces $\rho^{(a,b)}$ and $\wt z$ are not 
perpendicular. 

We consider
an orthonormal basis $\{\wt K_1,\dots ,\wt K_{k-3}\}$ 
also on $K^{(a,b)}$ 
and
the
matrix of the Ricci operator $\wt r$ is computed 
with respect to the basis $\{\wt K_1,\dots ,\wt K_{k-3},\wt E_{\bold i},
\wt E_{\bold j},\wt E_{\bold k},\bold j,\bold k\}$, in the 
nilpotent case. 
Then in \cite{Sz4} we get that this matrix is of the form
   $$\wt r=\pmatrix
     \epsilon I_K &
     0 & 0 & 0\\
     0 &(\epsilon +E_{ll})I_{l}& 0 & 0\\ 
     0 & 0 &
    (\epsilon 
+E_{LL})
I_{L} & 
     E_{{L}\wt z} \\
     0 & 0 &
     E_{\wt zL} &
    (\epsilon 
+E_{\wt z\wt z})
I_{\wt z}
     \endpmatrix,\tag 6.5$$
where $I_K$, $I_{\wt z}$ 
are unit matrices on
the spaces $K$, $\wt z$ 
and $I_{L}$ (resp. $I_{l}$) is unit matrix
 on the space ${L}$
 spanned by the vectors
$\bold j$ and $\bold k$ (resp. on
 the 1-dimensional space 
 $l$ 
spanned by $\bold i$).
 The entries $\epsilon\, ,\, E_{ll}\, ,\, E_{LL}\, ,\,
E_{\wt z\wt z}$ are very complicated expressions depending on
the function $\delta (|Z|)$ (cf. (4.20) in \cite{Sz4}). Yet it turns 
out that the determinant
   $$ det\pmatrix
     E_{LL}I_{L} &
     E_{L\wt z} \\
     E_{\wt zL} &
     E_{\wt z\wt z} I_{\wt z}
     \endpmatrix
\tag 6.6$$
is non-zero
 on an everywhere dense open set. This statement proves that on this set
the Ricci tensor has distinct eigenvalues on the
invariant subspaces $K$ and $\rho\oplus\wt z$. 

In the solvable case the Ricci curvature $\wt r_S$ 
 is computed
with respect to the basis $\{\wt K_1,\dots ,\wt K_{k-3},\wt E_{\bold i},
\wt E_{\bold j},\wt E_{\bold k},\bold j,\bold k,\bold t\}$. 
Then we get
$$
   \wt r_S=\pmatrix
     \sigma I_K &
     0 & 0 & 0 & 0\\
     0 &(\sigma +S_{ll})& 0 & 0 & S_{l\bold t}
\\ 
     0 & 0 &
    (\sigma 
+S_{LL})
I_{L} & t^{1\over2} 
     S_{L\wt z} &0 \\
     0 & 0 &t^{1\over2}
     S_{\wt zL} &
    (\sigma 
+S_{\wt z\wt z})
I_{\wt z} &0 \\
0&S_{\bold tl}&0&0& \sigma+S_{\bold t\bold t}
     \endpmatrix,\tag 6.7
$$

By the very same proof given in the nilpotent case we get that the
eigenvalue $\sigma$ on the eigenspace $K$ is different from the other
eigenvalues on an everywhere dense open set of a geodesic
sphere and therefore this 
distribution is invariant by the actions of the isometries. 

The {\it First Non-Isometry Proof} can be established as follows

 For the vector fields $U$ and $V$ tangent to
$\rho^{(a,b)}\oplus \wt z$ in the nilpotent case, resp. tangent to 
$\rho^{(a,b)}\oplus \wt z\oplus\wt T$ in the solvable case,  
 let $L(U,V)$ be the orthogonal
projection of $[U,V]$ onto $K^{(a,b)}$. The $L$ is obviously a tensor
field of type (2,1) on $\pa D$ and
it is invariant with respect to the action of isometries
on the space. By (2.4), this $L$ vanishes exactly at the 
points of the form 
$(X^{(a)},Z)$ or $(X^{(b)},Z)$, i. e. the induced metrics
on the hypersurface $\pa D$ of the spaces
$\bold H^{(a,b)}_3$ and
$\bold H^{(a',b')}_3$ (resp. on their solvable extensions)
can not be isometric unless $(a,b)=(a',b')$ up to an order.

This Proof also proves that the isometries on the geodesic spheres of
the spaces
$S\bold H^{(a,b)}_3$ with $ab\not =0$ are locally inhomogeneous and
thus the Non-Isometry Proof is completely established for the
Striking Examples.

In \cite{Sz4} we describe two other proofs.
In the {\it second Non-Isometry Proof} we trace back the problem to
the ambient space by showing that the considered surfaces are isometric
if and only if the ambient spaces are isometric. 

First we describe a geometric interpretation of the distribution 
$K^{(a,b)}$.
On the space $H^{(a,0)}_3$ (resp. on
 $H^{(0,b)}_3$)
 the distribution $\rho^{(a,0)}$ 
(resp. $\rho^{(0,b)}$) 
is integrable
and the 3-dimensional integral manifolds are the fibers of
 a principal fiber
bundle with the structure group $SO(3)$. 
This fibration is nothing but
the quaternionic Hopf fibration and the factor space is the
2-point homogeneous quaternionic projective space \cite {Be}.

If $a\, >\,1$
(resp. $b\, >\,1$),
the distribution $K^{(a,0)}$ (resp. $K^{(0,b)}$)
 is an irreducible connection
on this bundle with an irreducible curvature form $\omega (X,Y)=
[X,Y]_\rho ;X,Y\in K$.
This proves that 
$[K^{(a,b)},K^{(a,b)}]_\rho =\rho^{(a,b)}$. 
Therefore the $K$ generates the whole tangent
space on a sphere $S_X(Z)$ (resp. on $S_X(Z,t)$) by Lie brackets.
Since the $K^{(a,b)}$ is invariant by the actions of isometries,
the spheres $S_X$ must be mapped to such spheres by isometries.
By (1.11) (resp. by the corresponding version on the solvable extension)
we get that this action between these spheres is nothing but  
restriction of an orthogonal transformation, defined
between the ambient X-spaces. These
orthogonal transformations define an extension of the considered
isometry onto the ambient space, which turns out to be an isometry
between the ambient spaces.

In the third {\it Non-Isometry Proof} 
we determine the isometries on the considered
hypersurfaces. It turns out that in case of $ab\not =0$,
 the group of isometries
is the non-transitive group
$O(\bold H^a)\times O(\bold H^b)\times SO(3)$ (semi direct product),
 in both of the nilpotent
and the solvable cases, 
where $O(\bold H^{c})$ is the 
quaternionic orthogonal group acting on $\bold H^{c}$. These groups are
obviously different from the transitive isometry groups acting on the
corresponding 2-point homogeneous space. 
By summing up we have
\proclaim {Theorem 6.2}  The metrics on the 
same radius
geodesic spheres  
are iso\-spectral on a discrete family,
$S\bold H^{(a,b)}_3$, 
 of solvable $\si$-equivalent spaces.
 These spheres have
different local geometries unless $(a,b)=(a',b')$ up to an order.

 The geodesic spheres on a space 
$S\bold H^{(a+b,0)}_3$
 are homogeneous, while the spheres on the other spaces are locally
inhomogeneous. This demonstrate the fact:
"One can not hear the local homogeneity
property even on the most simple closed manifolds, namely on
 spheres." 
\endproclaim

\enddocument


\Refs\nofrills{\titlebf References}
\tenpoint
\widestnumber\key{\tenrm CDKR]}

\ref\key Be
\by\smc A. Besse
\book\rm Manifolds all of whose geodesics are closed
\publ Springer-Verlag Berlin Heidelberg New York
\yr 1978
\endref

\ref\key BT
\by\smc R. Brooks, R. Tse
\paper\rm Isospectral surfaces of small genus
\jour Nagoya Math. J.
\vol\rm 107
\yr 1987
\pages 13--24
\endref

\ref\key Bu
\by\smc P. Buser
\paper\rm Isospectral Riemann surfaces,
\jour Ann. Inst. Fourier (Grenoble)
\vol\rm 36
\yr 1986
\pages 167--192
\endref

\ref\key CH
\by\smc R. Courant, D. Hilbert
\book\rm Methods of mathematical physics,
Vol. I
\publ Interscience
\yr 1953
\endref

\ref\key Cha
\by\smc I. Chavel
\book\rm Eigenvalues in Riemannian geometry
\publ Academic Press
\yr 1984
\endref

\ref\key CDKR
\by\smc M. Cowling, A. H. Dooley, A. Koranyi, F. Ricci
\paper\rm $H$-type groups and Iwasawa decompositions
\jour Adv. in Math.
\vol\rm 87
\yr 1991
\pages 1--41
\endref

\ref\key CS
\by\smc C. Croke, V. Sharafutdinov
\paper\rm Spectral rigidity of a compact negatively curved manifold
\jour Topology
\vol\rm 37
\yr 1998
\pages 1265--1273
\endref

\ref\key D
\by\smc E. Damek
\paper\rm Geometry of a semidirect extension of a Heisenberg
type nilpotent group
\jour Colloq. Math.
\vol\rm 53
\yr 1987
\pages 255--268
\endref

\ref\key DR
\by\smc E. Damek, F. Ricci
\paper\rm A class of non-symmetric harmonic manifolds
\jour Bull. Amer. Math. Soc.
\vol\rm 27
\yr 1992
\pages 139--142
\endref


\ref\key DG
\by\smc D. DeTruck, C. Gordon
\paper\rm Isospectral deformations II.
\jour Com. Pure Appl. Math.
\vol\rm 42
\yr 1989
\pages 1067--1095
\endref

\ref\key E
\by\smc P. Eberlein
\paper\rm Geometry of two-step nilpotent groups with a left invariant metric
\jour Ann. Sci. \'Ecole Norm. Sup.
\vol\rm 27
\yr 1994
\pages 611--660
\endref

\ref\key G1
\by\smc C. Gordon
\paper\rm Isospectral closed Riemannian manifolds which are not
locally isometric
\jour J. Diff. Geom.
\vol\rm 37
\yr 1993
\pages 639--649
\endref

\ref\key G2
\by\smc C. Gordon
\paper\rm Isospectral closed manifolds which are not locally
isometric, Part II
\inbook Contemporary Mathematics: Geometry of the spectrum
\eds R. Brooks, C. Gordon, P. Perry
\publaddr vol. 173, AMS, 1994, 121--131
\endref

\ref\key GGSW
\by\smc C. Gordon, R. Gornet, D. Schueth, D. Webb, E. Wilson
\paper\rm Isospectral deformations of closed Riemannian manifolds with
different scalar curvature
\jour Ann. Inst. Four.,Grenoble
\vol\rm 48
\yr 1998
\pages 593--607
\endref

\ref\key GS
\by\smc C. Gordon, Z.I. Szab\'o
\paper\rm Isospectral deformations on negatively curved manifolds with 
boundaries
\jour Preprint
\endref

\ref\key GW1
\by\smc C. Gordon, E. N. Wilson
\paper\rm Isospectral deformations of compact solvmanifolds
\jour J. Diff. Geom.
\vol\rm 19
\yr 1984
\pages 241--256
\endref

\ref\key GW2
\by\smc C. Gordon, E. N. Wilson
\paper\rm The spectrum of the Laplacian on Riemannian Heisenberg
manifolds
\jour Michigan Math. J.
\vol\rm 33
\yr 1986
\pages 253--271
\endref

\ref\key GW3
\by\smc C. Gordon, E. N. Wilson
\paper\rm Continuous families of isospectral Riemannian metrics
which are not locally isometric
\jour J. Diff. Geom.
\vol\rm 47
\yr 1997
\pages 504--529
\endref

\ref\key GWW
\by\smc C. Gordon, D. Webb, S. Wolpert
\paper\rm One can not hear the shape of a drum
\jour Bull. Amer. Math. Soc.
\vol\rm 27, no. 1
\yr 1992
\pages 134--138
\endref

\ref\key I
\by\smc A. Ikeda
\paper\rm On lens spaces which are isospectral but not isometric
\jour Ann. Sci. Ecole Norm. Sup.
\vol\rm (4)13
\yr 1980
\pages 303--315
\endref

\ref\key K
\by\smc A. Kaplan
\paper\rm Riemannian nilmanifolds attached to Clifford modules
\jour Geom. Dedicata
\vol\rm 11
\yr 1981
\pages 127--136
\endref

\ref\key Ka
\by\smc R. Karidi
\paper\rm Ricci curvature and volume growth for leftinvariant Riemannian
metrics on nilpotent and some solvable Lie groups
\jour Geom. Dedicata
\vol\rm 46
\yr 1993
\pages 249--277
\endref

\ref\key M
\by\smc J. Milnor
\paper\rm Eigenvalues of the Laplace operator on certain manifolds
\jour Proc. Nat. Acad. Sci. U.S.A.
\yr 1964
\pages 542
\endref

\ref\key P
\by\smc H. Pesce
\paper\rm Calcul du spectre d'une nilvariete de rang deux et
applications
\jour Trans. Amer. Math. Soc.
\vol\rm 339
\yr 1993
\pages 433--461
\endref

\ref\key Rh
\by\smc C. Riehm
\paper\rm Explicit spin representation and Lie algebras of
Heisenberg type
\jour J. London Math. Soc.
\vol\rm 29
\yr 1984
\pages 49--62
\endref

\ref\key Sch
\by\smc D. Schueth
\paper\rm Continuous families of isospectral metrics on simply connected 
manifolds
\jour Ann.of Math.
\vol\rm 149
\yr 1999
\pages 287--308
\endref

\ref\key S
\by\smc T. Sunada
\paper\rm Riemannian coverings and isospectral manifolds
\jour Ann. of Math.
\vol\rm (2)121
\yr 1985
\pages 169--186
\endref

\ref\key Sz1
\by\smc Z. I. Szab\'o
\paper\rm Lichnerowicz conjecture on harmonic manifolds
\jour J. Diff. Geom.
\vol\rm 31
\yr 1990
\pages 1--28
\endref

\ref\key Sz2
\by\smc Z. I. Szab\'o
\paper\rm Spectral theory for operator families on Riemannian
Manifolds
\inbook Proc. of Symp. in Pure Math. vol. 54 (1993), Part 3
\yr 615--665
\endref

\ref\key Sz3
\by\smc Z. I. Szab\'o
\paper\rm Locally non-isometric yet super isospectral spaces
\jour Geom. funct. anal. (GAFA)
\vol\rm 9
\yr 1999 (in preprint 1992)
\pages 185--214
\endref

\ref\key Sz4
\by\smc Z. I. Szab\'o
\paper\rm Cornucopia of isospectral pairs of metrics on
balls and spheres with different local geometries
\jour Preprint
\yr 2000 
\endref

\ref\key V
\by\smc M.-F. Vigneras
\paper\rm Varietes Riemanniennes isospectrales et non\linebreak
isometrique
\jour Ann. of Math.
\vol\rm (2)112
\yr 1980
\pages 21--32
\endref

\ref\key W
\by\smc E. N. Wilson
\paper\rm Isometry groups on homogeneous nilmanifolds
\jour Geom. Dedicata
\vol\rm 12
\yr 1982 
\pages 337--346
\endref

\endRefs
\bye